\documentclass[12pt]{amsart}
\usepackage{a4wide}

\usepackage{enumerate}
\usepackage{booktabs, multirow}
\usepackage{amsmath,amsfonts,amssymb,amsthm}
\usepackage{xcolor}
\usepackage{pb-diagram}
\usepackage{array, hhline, url}

\newtheorem{theorem}{Theorem}[section]
\newtheorem{lemma}[theorem]{Lemma}
\newtheorem{proposition}[theorem]{Proposition}

\newtheorem{innercustomthm}{Theorem}
\newenvironment{customthm}[1]
{\renewcommand\theinnercustomthm{#1}\innercustomthm}
{\endinnercustomthm}

\theoremstyle{definition}
\newtheorem{remark}[theorem]{Remark}

\newcommand{\R}{\ensuremath{\mathbb{R}}}
\renewcommand{\H}{\ensuremath{\mathbb{H}}}
\newcommand{\g}[1]{\ensuremath{\mathfrak{#1}}}

\DeclareMathOperator{\tr}{tr}
\DeclareMathOperator{\Isom}{Isom}
\DeclareMathOperator{\sech}{sech}

\DeclareMathOperator{\Ad}{Ad}
\DeclareMathOperator{\ad}{ad}
\DeclareMathOperator{\id}{id}

\DeclareMathOperator{\Exp}{Exp}

\begin{document}
	
\title{Polar actions on homogeneous $3$-spaces}

\author[M.\ Dom\'{\i}nguez--V\'{a}zquez]{Miguel Dom\'{\i}nguez-V\'{a}zquez}
\author[T.\ A.\ Ferreira]{Tarcios A.\ Ferreira}
\author[T.\ Otero]{Tom\'as Otero}

\address{Department of Mathematics, Universidade de Santiago de Compostela, Spain 
\newline \indent CITMAga, 15782 Santiago de Compostela, Spain}
\email{miguel.dominguez@usc.es}
\address{Departamento de Matem\'atica, Universidade de Brasi\'ilia, Brazil}
\email{tarciosandrey@gmail.com}
\address{Mathematisches Institut, Universität Münster, Einsteinstr.~62, 48149 Münster, Germany}
\email{tomas.otero@uni-muenster.de}

\thanks{The first and third authors have been supported by grant PID2022-138988NB-I00 funded by MICIU/AEI/10.13039/501100011033 and by ERDF, EU, and by project ED431C 2023/31 (Xunta de Galicia, Spain). The second author has been supported by CAPES Print Process N° 88887.890396/2023-00. The third author acknowledges support from from the Deutsche Forschungsgemeinschaft (DFG, German Research Foundation) Project ID~427320536--SFB 1442, and Germany’s Excellence Strategy EXC~2044--390685587, Mathematics Münster: Dynamics-Geometry-Structure.}

\subjclass[2020]{Primary 53C30; Secondary 53C40}


\begin{abstract}
We classify polar isometric actions on simply connected $3$-dimensional Riemannian homogeneous spaces, up to orbit equivalence. In particular, we classify extrinsically homogeneous surfaces in such spaces and study the geometry of the orbit foliations of the corresponding cohomogeneity one actions.
\end{abstract}

\keywords{Polar action, homogeneous surface, cohomogeneity one action, homogeneous space, constant mean curvature}

\maketitle

 
\section{Introduction}\label{sec:intro}
A proper isometric action of a connected Lie group $H$ on a Riemannian manifold $M$ is said to be \emph{polar} if there exists a submanifold $\Sigma$ of $M$ which meets all the $H$-orbits, and every intersection of an $H$-orbit with $\Sigma$ is orthogonal. 
Such a $\Sigma$ is called a section of the polar action, and it must be a totally geodesic submanifold of $M$. The minimal codimension of the $H$-orbits, which coincides with the dimension of $\Sigma$, is called the cohomogeneity of the action.

Polar actions and their sections were proposed by Palais and Terng~\cite{PT} as generalizations of polar coordinates and canonical forms, respectively, to arbitrary Riemannian manifolds with a large isometry group. As they already observed, polar actions are interesting from different viewpoints, such as in invariant theory, Riemannian submanifold geometry and calculus of variations. But they turned out to be crucial in other problems, such as the construction of geometric structures on manifolds and the study of positively curved spaces, as in~\cite{FH, FGT}. This is especially true for the particular case of actions with one-dimensional orbit space (which are necessarily polar), commonly referred to as \emph{cohomogeneity one actions}. On the other hand, the consideration of cohomogeneity two polar actions allows for the construction of Riemannian submanifolds with specific properties, via the so-called equivariant method, as in~\cite{HH:inventiones} or~\cite{AR:acta}.

The motivation for this article comes, ultimately, from the viewpoint of Riemannian submanifold geometry. The families of orbits of polar actions on Riemannian manifolds enjoy several interesting properties, making their study worthwhile. For instance, they constitute singular Riemannian foliations whose regular orbits have parallel mean curvature. In particular, the regular orbits of cohomogeneity one actions are precisely the \emph{(extrinsically) homogeneous hypersurfaces}, that is, the hypersurfaces such that for any two of their points there is an isometry of the ambient space mapping one point to the other. Therefore, homogeneous hypersurfaces have constant principal curvatures and, in particular, constant mean curvature. Moreover, the orbit foliation of a cohomogeneity one action is an isoparametric family of hypersurfaces, that is, a decomposition of the ambient space into parallel submanifolds, all of them of codimension one and with constant mean curvature, except for at most two singular orbits (which have higher codimension). 

Thus, the classification of polar actions on a given ambient space, up to orbit equivalence, constitutes a natural problem in submanifold geometry. We recall that two isometric actions on a Riemannian manifold $M$ are called orbit equivalent if there is an isometry of $M$ mapping the orbits of one action to the orbits of the other one. This classification problem is particularly relevant on those manifolds with a large isometry group, such as homogeneous spaces. In fact, it has been deeply investigated in the special case of symmetric spaces. We refer to~\cite[Chapters~12-13]{BCO} for details and references. For example, the classification in space forms is known, due to works of Dadok~\cite{Dadok} and Wu~\cite{Wu}. Actually, the classification of polar actions on irreducible symmetric spaces of compact type is known. However, there are still many open questions, such as the behavior of polar actions on reducible symmetric spaces or the classification in most irreducible spaces of noncompact type. Indeed, cohomogeneity one actions on symmetric spaces of noncompact type were only classified very recently~\cite{SoSa}. The problem in the more general setting of homogeneous spaces is essentially untackled; see~\cite{GKRV} for a recent contribution in the framework of certain Stiefel manifolds.

In this article, we classify polar actions on simply connected Riemannian homogeneous $3$-spaces up to orbit equivalence. Note that, since the ambient dimension is $3$, the only nontrivial and nontransitive polar actions have cohomogeneity one or two. Moreover, our description of cohomogeneity one actions will allow us to describe the geometry of the isoparametric families of homogeneous surfaces of these ambient spaces. 

In order to state our results, we will recall below some known facts about $3$-dimensional homogeneous spaces. Thus, let $M$ be a simply connected homogeneous Riemannian manifold of dimension $3$.
Then, it is well known that its isometry group $\Isom(M)$ must have dimension $3$, $4$ or $6$.

If $\dim(\Isom(M))=6$, then $M$ is a space form.
It essentially follows from Dadok's work \cite{Dadok} that any polar action on a Euclidean space $\R^n$ is orbit equivalent to the product of the isotropy representation of a symmetric space times a translational part, while polar actions on round spheres $\mathbb{S}^n(\kappa)$ are restrictions of polar representations on $\R^{n+1}$ to the sphere $\mathbb{S}^n(\kappa)$ centered at the origin. Hereafter, by $\mathbb{S}^n(\kappa)$ (resp.\ $\mathbb{H}^n(\kappa)$) we mean the sphere (resp.\ real hyperbolic space) of constant sectional curvature $\kappa$. The classification of polar actions on real hyperbolic spaces $\H^n(\kappa)$ is due to Wu~\cite{Wu}, and can be reduced to the classification on round spheres $\mathbb{S}^k$, $k<n$. See Remarks~\ref{rem:C1_Ekt} and~\ref{rem:polar_space_forms} for further details.

Simply connected homogeneous $3$-spaces with $\dim(\Isom(M))=4$ are isometric to bundles $\mathbb{E}(\kappa,\tau)$ over a complete, simply-connected surface of constant curvature $\kappa$ and bundle curvature $\tau$ with $\kappa-4\tau^2\neq 0$. The classification of cohomogeneity one actions on $\mathbb{E}(\kappa,\tau)$-spaces follows from the classification of isoparametric surfaces due to Manzano and the first author~\cite{DVMa:annali}. However, as far as the authors know, the classification of cohomogeneity two polar actions on these spaces has not been obtained yet. Such classification is part of the content of Theorem~B below, and will be addressed in Section~\ref{sec:cohom_2}.

Lastly, if $\dim(\Isom(M))=3$, then $M$ is isometric to a metric Lie group, that is, a Lie group $G$ endowed with a left-invariant metric. In this paper, we obtain the explicit classification of polar actions on these spaces, up to orbit equivalence. It is important to mention that, in the case of cohomogeneity one actions, the corresponding classification in these spaces is intimately related to the classification of $2$-dimensional subgroups of the corresponding groups $G$, which was addressed by Meeks and Pérez~\cite{MePe:cont}. However, they do not develop a systematic study of the orbit equivalence of the actions. Thus, in Section~\ref{sec:c1} we provide an alternative approach that yields the explicit moduli space of actions up to orbit equivalence. The result is stated in Theorem~A and Table~\ref{table C1} below, where the Lie algebra $\g{g}$ of $G$ and the Lie algebras $\g{h}$ of the subgroups $H$ of $G$ acting with cohomogeneity one are given in terms of an orthonormal basis $\{E_1,E_2,E_3\}$ of $\g{g}$; note also that we denote by $e$ the identity element of $G$. The Lie algebra structure of each $\g{g}$ is summarized in the second column of Table~\ref{table C1} (the remaining brackets are zero or can be deduced from the information given). We also include information about the extrinsic geometry of the associated homogeneous surfaces in the last column. This information will require a specific investigation of the family of orbits of each one of the actions in terms of parallel normal displacement, see Section~\ref{sec:geometry}.

We can now state our main results. The first one, along with previous known results for space forms and $\mathbb{E}(\kappa,\tau)$-spaces, gives the classification of cohomogeneity one actions (or equivalently, of homogeneous surfaces) of simply connected homogeneous $3$-spaces.

\begin{customthm}{A}\label{th:c1}
	An isometric action of a connected Lie group on a $3$-dimensional simply connected metric Lie group $G$ 
	with $\dim\left(\Isom(G)\right)=3$ is of cohomogeneity one if and only if it is orbit equivalent to the action of the connected Lie subgroup $H$ of $G$ with one (and only one) of the Lie algebras $\g{h}$ in Table~\ref{table C1}.
\end{customthm}

\renewcommand{\arraystretch}{1.2}
\begin{table}[h]
	\centering
	\begin{tabular}{l>{\footnotesize}l>{\footnotesize}l>{\footnotesize}l}
		\toprule
		$G$ & {\normalsize $\g{g}$} & {\normalsize $\g{h}$} & {\normalsize Orbits}  \\ \midrule[0.08em]
		\multirow{3}{*}{$\widetilde{\mathsf{E}}_2$} & $[E_2,E_3]=\lambda_1 E_1$ & \multirow{3}{*}{$\operatorname{span}\{E_1,E_2\}$} & \multirow{3}{*}{\parbox{4.2cm}{All orbits are minimal, but not totally geodesic.}}\\
		& $[E_3,E_1]=\lambda_2 E_2$ & & \\ 
		& $\lambda_1>\lambda_2>0$ & & \\ \hline
		\multirow{4}{*}{$\widetilde{\mathsf{SL}_2}(\R)$} & $[E_2,E_3]=\lambda_1 E_1$ & \multirow{4}{*}{$\operatorname{span}\{\sqrt{\lambda_1}E_1+\sqrt{-\lambda_3} E_3,E_2\}$} & \multirow{4}{*}{\parbox{4.2cm}{All orbits are minimal.\\No orbit is totally geodesic, except $H\cdot e$ precisely when $\lambda_2=\lambda_1+\lambda_3$. }}\\
		& $[E_3,E_1]=\lambda_2 E_2$ & & \\
		& $[E_1,E_2]=\lambda_3 E_3$ & & \\  
		& $\lambda_1>\lambda_2>0>\lambda_3$ & & \\ \hline
		\multirow{6}{*}{$\mathsf{Sol}_3$} &  & \multirow{3}{*}{$\operatorname{span}\{E_1,E_3\}$} & \multirow{3}{*}{\parbox{4.2cm}{All orbits are minimal, but not totally geodesic.}} \\
		& \multirow{2}{*}{$[E_2,E_3]=\lambda_1 E_1$}& & \\
		& \multirow{2}{*}{$[E_1,E_2]=\lambda_3 E_3$} & & \\ \cline{3-4}
		& & \multirow{3}{*}{$\operatorname{span}\{\sqrt{\lambda_1}E_1+\sqrt{-\lambda_3} E_3,E_2\}$} &  \multirow{3}{*}{{\parbox{4.2cm}{All orbits are minimal.\\ No orbit is totally geodesic, except $H\cdot e$ when $\lambda_3=-\lambda_1$.}}} \\
		& \multirow{1}{*}{$\lambda_1>0>\lambda_3$} & &\\
		& & & \\ \hline
		\multirow{10}{*}{$\R^2\rtimes\R$} &    & \multirow{2}{*}{$\operatorname{span}\{E_1,E_2\}$}    & \multirow{2}{4.2 cm}{All orbits have constant mean curvature $1$.}    \\
		&\multirow{6}{*}{$[E_2,E_3]=(1-\alpha)(\beta E_1-E_2)$}    						& & \\  \cline{3-4}
		& 							& $\operatorname{span}\{E_1+\frac{\alpha}{(1-\alpha)\beta} E_2,E_3\}$,    & \multirow{10}{*}{\parbox{4.2cm}{Each orbit has constant mean curvature in $(-1,1)$, and the map sending each orbit to its mean curvature is a bijection from the orbit space to $(-1,1)$.
			\\ 
		The orbit $H\cdot e$ is minimal, and it is totally geodesic if and only if $\beta=0$.} }   \\ 
		&\multirow{4}{*}{$[E_3,E_1]=(1+\alpha) (E_1+\beta E_2)$}   & when $\det L=1$    &    \\  \cline{3-3}
		& 							& $\operatorname{span}\{E_1,E_3\}$, when $\beta=0$   &    \\ \cline{3-3}
		&\multirow{3}{*}{$\alpha,\beta\geq 0$, $\alpha\neq 0,1$}   &  $\operatorname{span}\{E_2,E_3\}$, when $\beta=0$  &  \\ \cline{3-3}
		& 	\multirow{3}{*}{$\det L =(1-\alpha^2)(1+\beta^2)$}						& \multirow{2}{*}{$\operatorname{span}\{E_1+ \frac{\alpha+\sqrt{1-\det L}}{(1-\alpha)\beta} E_2,E_3\}$,}    &    \\  &&& \\
		& 	& when $\beta\neq 0$ and $\det L<1$ & \\ \cline{3-3}
		& 							& \multirow{2}{*}{$\operatorname{span}\{E_1+ \frac{\alpha-\sqrt{1-\det L}}{(1-\alpha)\beta} E_2,E_3\}$,}   &   \\ &&&\\
		&     						& when $\beta\neq 0$ and $\det L<1$ &  \\ 
		\bottomrule
	\end{tabular}	
	\bigskip
	\caption{Cohomogeneity one actions on $3$-dimensional metric Lie groups with $\dim(\Isom(G))=3$, up to orbit equivalence.}
	\label{table C1}
\end{table}

Our second theorem addresses the classification of polar actions of cohomogeneity two on any simply connected homogeneous $3$-space.

\begin{customthm}{B}\label{th:c2}
	An isometric action of a connected Lie group on a $3$-dimensional simply connected homogeneous space $M$ is polar and of cohomogeneity two if and only if one of the following holds:
	\begin{itemize}
		\item $M$ is a Riemannian symmetric space and the action is orbit equivalent to precisely one of the well-known $H$-actions in Table~\ref{table C2:symm spaces}.
		\item $M$ is isometric to a metric Lie group $G$ with $\dim(\Isom(G))=3$ and the action is orbit  equivalent to the action of the connected Lie subgroup $H$ of $G$ with precisely one of the Lie subalgebras $\g{h}$ in Table \ref{table C2:groups}.
	\end{itemize}
\end{customthm}

\begin{table}[h]
	\centering
	\begin{tabular}{l@{\qquad}l@{\qquad}l}
		\toprule
		$G$ & $\g{g}$ & $\g{h}$ \\ \midrule[0.08em]
		\multirow{3}{*}{$\mathsf{Sol}_3$} & $[E_2,E_3]=\lambda_1 E_1$ & \multirow{3}{*}{$\operatorname{span}\{E_1-E_3\}$}\\
		& $[E_1,E_2]=-\lambda_1 E_3$ & \\ 
		& $\lambda_1>0$ &\\ \hline
		\multirow{3}{*}{$\R^2\rtimes \R$} & $[E_2,E_3]=(\alpha-1) E_2$ & \multirow{2}{*}{$\operatorname{span}\{E_1\}$}\\ 
		& $[E_3,E_1]=(1+\alpha) E_1$ & \\
		&{$\alpha\neq 0,1$} & \multirow{1}{*}{$\operatorname{span}\{E_2\}$}  \\
		 \bottomrule
	\end{tabular}
	\bigskip
	\caption{Cohomogeneity two polar actions on $3$-dimensional metric Lie groups with $\dim(\Isom(G))=3$, up to orbit equivalence.}
	\label{table C2:groups}
\end{table}

This article is organized as follows. Section~\ref{sec:structure} is devoted to review the classification of simply connected $3$-dimensional homogeneous spaces, with particular emphasis on the subclass of metric Lie groups. In Section~\ref{sec:c1} we derive the classification of cohomogeneity one polar actions up to orbit equivalence on metric Lie groups, thus proving Theorem~\ref{th:c1}. In Section~\ref{sec:cohom_2} we prove Theorem~\ref{th:c2}, which contains the classification of cohomogeneity two polar actions. Finally, Section~\ref{sec:geometry} is devoted to the study of the extrinsic geometry of the families of parallel homogeneous surfaces induced by each one of the cohomogeneity one actions classified in Section~\ref{sec:c1}.
	

\section{Riemannian homogeneous $3$-manifolds and metric Lie groups}\label{sec:structure}

In this section we recall the classification of simply connected Riemannian homogeneous $3$-manifolds. For further information we refer to~\cite{MaSo:mathZ}, \cite{MePe:cont} and~\cite{Mi:adv}.

A Riemannian metric on a Lie group $G$ is said to be left-invariant if the left multiplication by $g$, $L_g$, is an isometry of $G$ for all $g\in G$. Therefore, a left-invariant metric on $G$ is determined by the choice of an inner product $\langle\cdot,\cdot\rangle$ on the Lie algebra $\g{g}$. We will denote by $\langle\cdot,\cdot\rangle$ both the inner product on $\g{g}$ and the associated Riemannian metric on $G$.
Lie groups equipped with a left-invariant metric are often called metric Lie groups, and provide perhaps the simplest examples of Riemannian homogeneous spaces.

In the special case dimension $3$, it turns out that  metric Lie groups are close to exhausting all examples of Riemannian homogeneous manifolds:
any simply connected homogeneous $3$-manifold different from $\mathbb{S}^2(\kappa)\times\R$ is isometric to a $3$-dimensional metric Lie group (cf. \cite[Theorem 2.4]{MePe:cont}).
It should be noted, however, that non-isomorphic Lie groups can give rise to isometric Riemannian homogeneous spaces. We also note that, among simply connected homogeneous $3$-manifolds, the subclass of Riemannian symmetric spaces is constituted by the space forms $\R^3$, $\mathbb{S}^3$ and $\mathbb{H}^3$ and the product spaces $\mathbb{S}^2\times\R$ and $\mathbb{H}^2\times\R$, where the latter correspond to $\mathbb{E}(\kappa,\tau)$-spaces with $\tau=0$ and $\kappa\neq 0$.

We will summarize below Milnor's description of $3$-dimensional Lie algebras \cite{Mi:adv}.
For the sake of convenience, one distinguishes two cases depending on the unimodularity of the group:
a Lie group $G$ is said to be unimodular if its left and right-invariant Haar measures coincide, and is called non-unimodular otherwise.
For a connected Lie group $G$, being unimodular is equivalent to the condition that the adjoint transformation $\ad_X\colon\g{g}\rightarrow\g{g}$, $\ad_X(\cdot)=[X,\,\cdot\,]$, has trace $0$ for every $X\in\g{g}$.
Thus, a Lie algebra $\g{g}$ satisfying $\tr(\ad_X)=0$ for all $X\in \g{g}$ is called a unimodular Lie algebra. 

\subsection{Unimodular Lie groups of dimension $3$}\label{subsec:unimodular}
Let $\g{g}$ be a $3$-dimensional unimodular Lie algebra with inner product $\langle\cdot,\cdot\rangle$.
Given an orientation in $\g{g}$, both the corresponding cross product $\wedge$ and Lie bracket $[\cdot,\cdot]$ are skew-symmetric bilinear operators on $\g{g}$. 
Thus, one can define a linear map $L\colon\g{g}\rightarrow\g{g}$ by the formula
$
[X,Y]=L(X\wedge Y)$,
which turns out to be self-adjoint by the unimodularity of $\g{g}$.
Therefore, there exists an orthonormal basis $\{E_1,E_2,E_3\}$ of $\g{g}$ such that
\begin{equation}\label{unimodular brackets}
	[E_1,E_2]=\lambda_3 E_3,\qquad [E_2,E_3]=\lambda_1 E_1,\qquad [E_3,E_1]=\lambda_2 E_2,
\end{equation}
for some real constants $\lambda_1,\lambda_2,\lambda_3$. These constants ultimately determine the geometry and group structure of the unique simply connected metric Lie group $G$ with Lie algebra $\g{g}$ and left-invariant metric determined by $\langle\cdot,\cdot\rangle$.

Note that changing the sign of all the $\lambda_i$ corresponds to a change in the orientation in $\g{g}$, but the metric structure does not change.
Similarly, multiplying all the $\lambda_i$ by a positive number corresponds to rescaling the metric in $G$, but the underlying group structure remains the same.
From now on, we will assume without loss of generality that $\lambda_1\geq\lambda_2\geq\lambda_3$ and at most one $\lambda_i<0$.
A list of the corresponding Lie groups is given in Table \ref{tab:unimgroups}.
\begin{table}[h]
	\centering
		\begin{tabular}{l@{\qquad}c@{\qquad}c@{\qquad}c@{\qquad}c@{\qquad}c@{\qquad}c}
		\toprule
		Lie group&  $\mathsf{SU}_2$ & $\widetilde{\mathsf{E}}_2$ & 
		$\widetilde{\mathsf{SL}_2}(\mathbb{R})$& $\mathsf{Sol}_3$& $\mathsf{Nil}_3$& $\R^3$
		\\		\midrule
		Signs of $\lambda_1,\lambda_2,\lambda_3$ & 
		$+,+,+$ & 
		$+,+,0$ &
		$+,+,-$ & 
		$+,0,-$ &
		$+,0,0$ & 
		$0,0,0$\\
		\bottomrule
	\end{tabular}
	\bigskip
	\caption{Three-dimensional simply-connected unimodular Lie groups in terms of the structure constants.} \label{tab:unimgroups}
\end{table} 

Let $\mu_i=\frac{1}{2}(\lambda_1+\lambda_2+\lambda_3)-\lambda_i$, for each $i=1,2,3$. The Koszul formula allows us to obtain the Levi-Civita connection ${\nabla}$ of $G$, which is given by
\begin{equation}\label{unimodular connection}
	\begin{aligned}
		{\nabla}_{E_1}E_1&=0, & {\nabla}_{E_1}E_2&=\mu_1 E_3, & {\nabla}_{E_1}E_3&=-\mu_1 E_2,  \\
		{\nabla}_{E_2}E_1&=-\mu_2 E_3, & {\nabla}_{E_2}E_2&=0, & {\nabla}_{E_2}E_3&=\mu_2 E_1, \\
		{\nabla}_{E_3}E_1&=\mu_3 E_2, & {\nabla}_{E_3}E_2&=-\mu_3 E_1, & {\nabla}_{E_3}E_3&=0.  
	\end{aligned}
\end{equation}

\begin{remark}\label{remark isometry group unimodular}
	If two of the structure constants are equal, say $\lambda_1=\lambda_2\neq\lambda_3$, and $\lambda_3\neq 0$, one can see that $(G,\langle\cdot,\cdot\rangle)$ is isometric to $\mathbb{E}(\lambda_1\lambda_3,\frac{\lambda_3}{2})$, and $\Isom(G)$ has dimension $4$.
	If $\lambda_1=\lambda_2\neq\lambda_3=0$ then $G$ is isomorphic to the universal cover $\widetilde{\mathsf{E}}_2$ of the Euclidean group of the plane (which in turn is isometric to the Euclidean space $\R^3$). If $\lambda_1=\lambda_2=\lambda_3$, then $G$ is isomorphic to $\R^3$ or $\mathsf{SU}_2$, and it has nonnegative constant sectional curvature and a $6$-dimensional isometry group.
	When the three structure constants are different, then $\Isom(G)$ has dimension $3$.
	A description of the full isometry group of each unimodular $G$ can be found in \cite{HaLee}.
\end{remark}

\subsection{Non-unimodular Lie groups of dimension $3$}\label{subsec:nonunim}
Let $\g{g}$ be a non-unimodular $3$-dimensional Lie algebra.
Then it is a semi-direct product $\g{g}=\R^2\rtimes\R$ and, up to rescaling the metric, there exists an orthonormal basis $\{E_1,E_2,E_3\}$ of $\g{g}$ such that $\R^2=\operatorname{span}\{E_1,E_2\}$, and $\ad_{E_3}$ is given by the matrix
\[
	L=\begin{pmatrix}
	(1+\alpha) & -(1-\alpha)\beta\\
	(1+\alpha)\beta & (1-\alpha)
\end{pmatrix}
\]
in the basis $\{E_1,E_2\}$, for some constants $\alpha,\beta\geq 0$.
That is, the Lie bracket of $\g{g}$ is given~by
\begin{equation}\label{brackets nonunimodular}
	[E_1,E_2]=0,\quad [E_2,E_3]=(1-\alpha)\beta E_1+(\alpha-1)E_2,\quad [E_3,E_1]=(1+\alpha)E_1+(1+\alpha)\beta E_2.
\end{equation}

If $L\neq \operatorname{Id}$, the determinant of $L$, $\det L=(1-\alpha^2)(1+\beta^2)$, provides a complete isomorphism invariant for the Lie algebra $\g{g}$, that is, if $\det(L)=\det(M)$ and $L\neq\operatorname{Id}\neq M$, then $\R^2\rtimes_{L}\R\cong \R^2\rtimes_{M}\R$ as Lie algebras (although their metrics may differ).

Using the Koszul formula, we get
\begin{equation}\label{eq:LC_non_unimodular}
	\begin{aligned}
		 {\nabla}_{E_1}E_1&=(1+\alpha) E_3, &  {\nabla}_{E_1}E_2&=\alpha \beta  E_3, &  {\nabla}_{E_1}E_3&=-(1+\alpha)E_1-\alpha \beta  E_2,  \\
		 {\nabla}_{E_2}E_1&=\alpha \beta  E_3, &  {\nabla}_{E_2}E_2&=(1-\alpha) E_3, &  {\nabla}_{E_2}E_3&=-\alpha \beta  E_1-(1-\alpha) E_2, \\
		 {\nabla}_{E_3}E_1&=\beta  E_2, &  {\nabla}_{E_3}E_2&=-\beta  E_1, &  {\nabla}_{E_3}E_3&=0.  
	\end{aligned}
\end{equation}

\begin{remark}\label{remark isometry group nonunimodular}
	If $\alpha=0$, then $G$ is isometric to $\H^3(-1)$, whereas if $\alpha=1$, then $G$ is isometric to $\mathbb{E}(-4,\beta )$, whose isometry group has dimension $4$. If $\alpha\notin\{0,1\}$ one can show that $\dim\Isom(G)=3$ by looking at the eigenvalues of the Ricci operator.
	A description of the full isometry group of $G$ can be found in \cite{CoRe}.
\end{remark}


\section{Cohomogeneity one actions on $3$-dimensional groups}\label{sec:c1}
The aim of this section is to classify cohomogeneity one actions on $3$-dimensional simply connected metric Lie groups, up to orbit equivalence. Along with the known result for space forms and $\mathbb{E}(\kappa,\tau)$-spaces (see Remarks~\ref{rem:C1_sf} and~\ref{rem:C1_Ekt} below), this completes the classification on all simply connected $3$-dimensional homogeneous spaces.

\begin{remark}[Cohomogeneity one actions on space forms]\label{rem:C1_sf}
	The classification of cohomogeneity one $H$-actions on simply connected $3$-dimensional space forms~$M$ is summarized in Table~\ref{table:C1_sf}. This result is  well known and can be obtained directly or as a consequence of classical results on isoparametric hypersurfaces by Segre and Cartan in the 30s (see for example~\cite[p.~84]{BCO} or~\cite[pp.~96-99]{CR:book}). Actually, isoparametric hypersurfaces in Euclidean and real hyperbolic spaces are homogeneous; this is no longer true for round spheres, but it still holds for $\mathbb{S}^3(\kappa)$. 
	
	In the case of $\mathbb{H}^3(\kappa)$ in Table~\ref{table:C1_sf}, we use certain notation coming from the Iwasawa decomposition of the simple Lie group $G=\mathsf{SO}_{1,3}^0$, which is the connected component of the identity of $\Isom(\mathbb{H}^3(\kappa))$. This result states that a semisimple Lie group $G$ is diffeomorphic to a product manifold $K\times A\times N$, see~\cite[pp.~340 and 343]{BCO}. In this case, $K\cong\mathsf{SO}_3$ is the isotropy at some basepoint $o$, $A\cong\R$ is certain $1$-dimensional subgroup of $G$, and $N$ is a $2$-dimensional abelian subgroup of $G$. The $K$-action on $\mathbb{H}^3(\kappa)$ fixes the basepoint $o$ and the other orbits are geodesic spheres around it, the $A$-action gives rise to a geodesic through the basepoint $o$ and equidistant curves to it, and the $N$-action produces a horosphere foliation whose common point at infinity is one of the limit points of the geodesic $A\cdot o$. The combination of the $A$-action with a certain rotational $1$-dimensional subgroup of $K$ (precisely, the connected component of the identity of the centralizer of $A$ in $K$), or with any $1$-dimensional subgroup of $N$ yields the second and third actions in Table~\ref{table:C1_sf} for $\mathbb{H}^3(\kappa)$, respectively.
	
	In all cases in Table~\ref{table:C1_sf}, all $H$-orbits are surfaces, except for at most two singular $H$-orbits. If there are  singular orbits, the $2$-dimensional orbits are tubes of different radii around each one of the singular orbits. By tube of radius $r$ we mean the subset of the ambient space obtained by traveling a fixed distance $r$ in all normal directions to a submanifold of codimension $\geq 2$. 
\end{remark}
	\begin{table}[h]
	\centering
	\begin{tabular}{lll}
		\toprule
		$M$ & $H$ & Orbits  \\ \midrule[0.08em]
		\multirow{3}{*}{$\R^3$} & $\R^2$ & Parallel affine planes \\ \cline{2-3}
		& $\mathsf{SO}_2\times\R$ & A straight line and coaxial cylinders around it
		\\ \cline{2-3}
		& $\mathsf{SO}_3$ & A fixed point and spheres around it \\  \hline
		\multirow{2}{*}{$\mathbb{S}^3(\kappa)$} & $\mathsf{SO}_3$ & Two  antipodal fixed points and tot.\ geodesic $2$-spheres around them \\ \cline{2-3}
		& $\mathsf{SO}_2\times \mathsf{SO}_2$ & Two maximal circles and tori around them
		\\ \hline
		\multirow{4}{*}{$\H^3(\kappa)$} & $K\cong\mathsf{SO}_3$ & A fixed point and geodesic spheres around it \\ \cline{2-3}
		& $\mathsf{SO}_2\times \R$ & A geodesic and tubes around it
		\\ \cline{2-3}
		& $\R\ltimes\R$ & A totally geodesic $\mathbb{H}^2(\kappa)$ and its equidistant surfaces \\ \cline{2-3}
		& $N\cong \R^2$ & Horosphere foliation \\ \bottomrule
	\end{tabular}	
	\bigskip	
	\caption{Cohomogeneity one actions on $3$-dimensional space forms up to orbit equivalence.}
	\label{table:C1_sf}
\end{table}

\begin{remark}[Cohomogeneity one actions on  $\mathbb{E}(\kappa,\tau)$-spaces]\label{rem:C1_Ekt} 
		The classification of cohomogeneity one actions in $\mathbb{E}(\kappa,\tau)$-spaces up to orbit equivalence is summarized in Table~\ref{table:C1_Ekt}. It can be obtained from the classification of isoparametric surfaces in these spaces~\cite{DVMa:annali}. As stated there, a complete hypersurface in one of these spaces is homogeneous if and only if it is isoparametric. 
		
		Recall that an $\mathbb{E}(\kappa,\tau)$-space is the total space of a fiber bundle $\pi\colon\mathbb{E}(\kappa,\tau)\to\mathbb{M}^2(\kappa)$ over a complete, simply-connected surface of constant curvature $\kappa$. The parameter $\tau$ codifies the curvature of the bundle. Thus $\mathbb{E}(\kappa,\tau)$ is a product $\mathbb{M}^2(\kappa)\times\R$ if and only if $\tau=0$. Moreover, one requires $\kappa-4\tau^2\neq 0$, so that $\dim(\Isom(\mathbb{E}(\kappa,\tau)))=4$. Actually, one of the Killing fields of $\mathbb{E}(\kappa,\tau)$ is tangent to the fibers of $\pi$, which then turns out to be a Killing (Riemannian) submersion.
		
		According to \cite{DVMa:annali}, the only examples of homogeneous surfaces of $\mathbb{E}(\kappa,\tau)$, $\kappa-4\tau^2\neq 0$, are: vertical cylinders over a complete curve of constant curvature in $\mathbb{M}^2(\kappa)$, a horizontal slice $\mathbb{M}^2(\kappa)\times\{t_0\}$, $t_0\in\R$, when $\tau=0$, or a so-called parabolic helicoid when $\kappa<0$. By vertical cylinder over a subset of the base $\mathbb{M}^2(\kappa)$ we understand the preimage of such a subset under $\pi$. Recall that a curve of constant curvature in $\mathbb{M}^2(\kappa)$ is a geodesic circle if $\kappa>0$; a circle or a straight line if $\kappa=0$; or a geodesic, an equidistant curve to a geodesic, a geodesic circle or a horocycle if $\kappa<0$. Finally, we refer to~\cite{DVMa:annali} for the explicit parameterization of parabolic helicoids (see also~\cite[Example~4.4]{Otero} for an alternative description in the case of $\mathbb{H}^2(\kappa)\times\R$).
	\end{remark}

	\begin{table}[h]
	\centering
	\begin{tabular}{lll}
		\toprule
		$M$ & $H$ & Orbits  \\ \midrule[0.08em]
		\multirow{2}{*}{$\mathbb{S}^2(\kappa)\times\R$} & $\mathsf{SO}_2\times \R$ & Two vertical lines and vertical cylinders around them \\ \cline{2-3}
		& $\mathsf{SO}_3$ & Parallel horizontal slices $\mathbb{S}^2(\kappa)\times\{t_0\}$
				\\ \hline
		\multirow{1}{*}{$\mathbb{S}^3_{\mathrm{Berger}}$}& $\mathsf{SO}_2\times\mathsf{SO}_2$& Two vertical circles and vertical tori around them		
		\\  \hline
		\multirow{5}{*}{\begin{tabular}{@{}l}$\mathbb{H}^2(\kappa)\times\R$\\$\widetilde{\mathsf{SL}_2}(\R)$\end{tabular}} & $\mathsf{SO}_2\times\R$ & Vertical line and vertical cylinders around it \\ \cline{2-3}
		& $\R\times\R$ & Vertical cylinders over a geodesic  and over its parallel curves
		\\ \cline{2-3}
		& $\R\times\R$ & Vertical cylinders over a horocycle foliation
		\\ \cline{2-3}
		& $\mathsf{SL}_2(\R)$ & Parallel horizontal slices $\mathbb{H}^2(\kappa)\times\{t_0\}$ (only for $\mathbb{H}^2(\kappa)\times\R$)
		\\ \cline{2-3}
		& $\R\ltimes\R$ & Parallel parabolic helicoids 
		\\ \hline
		\multirow{2}{*}{$\mathsf{Nil}_3$} & $\mathsf{SO}_2\times\R$ & Vertical geodesic and tubes around it \\ \cline{2-3}
		& $\R\times\R$ & Vertical cylinders over a foliation by parallel lines
		\\ \bottomrule
	\end{tabular}	
	\bigskip	
	\caption{Cohomogeneity one actions on $\mathbb{E}(\kappa,\tau)$-spaces of nonconstant curvature.}
	\label{table:C1_Ekt}
\end{table}

Let $G$ be a simply connected $3$-dimensional Lie group with a left-invariant metric $\langle\cdot,\cdot\rangle$ and Lie algebra $\g{g}$.
Since the action of $G$ on itself by left translations is a proper isometric free action, a subgroup $H$ of $G$ acts on $G$ with orbits of dimension $\dim H$.
For an arbitrary metric Lie group one can recover its isometry group as the group of products
\[
\Isom(G,\langle\cdot,\cdot\rangle)=G\cdot\Isom (G, \langle\cdot,\cdot\rangle)_e,
\]
where $G$ is identified with its left translations and $\Isom (G, \langle\cdot,\cdot\rangle)_e$ denotes the group of isometries of $G$ fixing the identity element.

Since cohomogeneity one actions on space forms and $\mathbb{E}(\kappa,\tau)$-spaces are well understood, from now on we will suppose that $\dim(\Isom(G))=3$.
Thus, $\Isom (G, \langle\cdot,\cdot\rangle)_e$ can be identified with a discrete (and hence finite) subgroup of $\mathsf{O}_3$ via the isotropy representation.
Let $\operatorname{Aut}(G)$ be the group of automorphisms of $G$.
It follows from \cite[Corollary 2.8]{HaLee} that $\Isom (G, \langle\cdot,\cdot\rangle)_e=\operatorname{Aut}(G)\cap\Isom(G, \langle\cdot,\cdot\rangle)$ is precisely the group of isometric automorphisms of $G$,
and so, by~\cite[p.~192]{HaLee}, we have
\[
\Isom(G,\langle\cdot,\cdot\rangle)=G\rtimes(\operatorname{Aut}(G)\cap\Isom(G, \langle\cdot,\cdot\rangle)).
\]
Thus, any  effective isometric action of a connected Lie group on $G$ is equivalent to the action by left translations of some connected Lie subgroup $H$ of $G$.
Moreover, two subgroups $H,\tilde{H}$ of $G$ give rise to orbit equivalent actions if and only if $H=C_{g}\circ \varphi (\tilde{H})$ for some $g\in G$ and $\varphi\in\operatorname{Aut}(G)\cap\Isom(G, \langle\cdot,\cdot\rangle)$, where $C_g$ denotes conjugation by $g$.
Altogether, we get:
\begin{proposition}\label{prop:1-1}
	Let $G$ be a simply connected metric Lie group with $\dim(\Isom(G))=3$. Then, there is a one-to-one correspondence between isometric actions of connected Lie subgroups on $G$ up to orbit equivalence and Lie subalgebras of $\g{g}$ up to conjugacy and isometric automorphisms.
\end{proposition}

\begin{remark}
	Let $H$ be a connected Lie subgroup of a $3$-dimensional simply connected Lie group $G$.
	If $G$ is solvable, then $H$ is closed in $G$ (see \cite[p.~670]{Che}).
	Otherwise, $G= \widetilde{\mathsf{SL}_2}(\R)$ or $G=\mathsf{SU}_2$, and every abelian subgroup of $G$ is one-dimensional and closed.
	It follows from \cite[Theorem~15]{Mal} that any connected subgroup $H$ of $G$ is closed. Therefore, the action of $H$ on $G$ by left translations is proper, which implies that its orbits are closed embedded submanifolds, see for example~\cite[pp.~66-67]{Michor}.
\end{remark}

In the following, we will classify codimension one Lie subalgebras of $\g{g}$ up to conjugacy and isometric automorphisms, for any $3$-dimensional metric Lie group $G$ with $\dim \Isom(G,\langle\cdot,\cdot\rangle)=3$.
For this, we will make use of the structure results of Section \ref{sec:structure} in order to give an explicit description of the corresponding Lie algebras $\g{g}$.

\subsection{Codimension one subalgebras of unimodular Lie algebras}
Recall that if $\g{g}$ is a unimodular $3$-dimensional Lie algebra, there exists an orthonormal basis $\{E_1,E_2,E_3\}$ of $\g{g}$ such that
\begin{equation*}\label{unimodular brackets 2}
[E_1,E_2]=\lambda_3 E_3,\qquad [E_2,E_3]=\lambda_1 E_1,\qquad [E_3,E_1]=\lambda_2 E_2,
\end{equation*}
with $\lambda_1\geq\lambda_2\geq\lambda_3$, and the corresponding Lie group structure can be recovered from the signs of the $\lambda_i$.
Let $\g{h}$ be a codimension one subalgebra of $\g{g}$, and write $\g{h}=\operatorname{span}\{A,B\}$ for some $A,B\in\g{g}$.
Write $A=a_1 E_1 +a_2 E_2+a_3 E_3$ and $B=b_1 E_1 +b_2 E_2+b_3 E_3$.
Then, $\g{h}$ is a subalgebra of $\g{g}$ if and only if
\begin{equation}\label{subalgebra condition unimodular}
0=\langle[A,B],A\wedge B\rangle=(a_2b_3-a_3b_2)^2\lambda_1+(a_3b_1-a_1b_3)^2\lambda_2+(a_1b_2-a_2b_1)^2\lambda_3.
\end{equation}

If $G=\R^3$ or $G=\mathsf{Nil}_3$, then any left-invariant metric on $G$ has an isometry group of dimension $\geq 4$.
Also, the only proper connected subgroups of $\mathsf{SU}_2$ are the trivial subgroup or $\mathsf{SO}_2$, which is one-dimensional.
Thus, we will exclude these cases from our study, and restrict to the Lie groups $\widetilde{\mathsf{E}}_2$, $\widetilde{\mathsf{SL}_2}(\R)$ and~$\mathsf{Sol}_3$.

\subsubsection{Codimension one subalgebras of $\g{e}_2$.}
For $\g{g}=\g{e}_2$, we have that $\lambda_1\geq\lambda_2>0$ and $\lambda_3=0$, so $\g{e}_2=\operatorname{span}\{E_1,E_2\}\rtimes\R E_3$.
Then $\g{h}$ is a subalgebra of $\g{e}_2$ if and only if
\[
(a_2b_3-a_3b_2)^2\lambda_1+(a_3b_1-a_1b_3)^2\lambda_2=0.
\]
Since both $\lambda_1\geq\lambda_2>0$, we have that $\g{h}$ is a subalgebra if and only if
\[
	a_2b_3-a_3b_2=0 \qquad\text{and}\qquad
	a_3b_1-a_1b_3=0.
\]
Thus, $\g{h}=\operatorname{span}\{E_1,E_2\}$ is the only codimension one subalgebra of $\g{e}_2$.

\subsubsection{Codimension one subalgebras of $\g{sl}_2(\R)$.}\label{subalgebras sl2r}
The Lie algebra $\g{sl}_2(\R)$ is a real simple Lie algebra. Since any $2$-dimensional Lie algebra is solvable, it follows that any $2$-dimensional subalgebra of $\g{sl}_2(\R)$ is maximal solvable and thus, by definition, a parabolic subalgebra of $\g{sl}_2(\R)$, see~\cite[\S13.2, p.~340]{BCO}. But since $\g{sl}_2(\R)$ has real rank equal to one, it has only one proper parabolic subalgebra, up to conjugacy \cite[p.~348]{BCO}.
Choosing $A=\sqrt{\lambda_1}\, E_1+\sqrt{-\lambda_3}\, E_3$ and $B=E_2$, with $\lambda_1>0>\lambda_3$, one immediately verifies~\eqref{subalgebra condition unimodular}:
\[
\langle[A,B],A\wedge B\rangle=\bigl(-\sqrt{-\lambda_3}\bigr)^2\lambda_1+\bigl(\sqrt{\lambda_1}\bigr)^2\lambda_3=-\lambda_1\lambda_3+\lambda_1\lambda_3=0.
\]
Hence $\g{h}=\operatorname{span}\{\sqrt{\lambda_1}\, E_1+\sqrt{-\lambda_3}E_3,E_2\}$ gives us a representative for this conjugacy class.

\subsubsection{Codimension one subalgebras of $\g{sol}_3$.}
For $\g{g}=\g{sol}_3$, we have $\lambda_1>\lambda_2=0>\lambda_3$. Hence, $\g{h}_0=\operatorname{span}\{E_1,E_3\}$ is a $2$-dimensional abelian subalgebra, and $\g{sol}_3=\g{h}_0\rtimes \R E_2$.

Suppose that $\g{h}\subset \g{sol}_3$ is some other $2$-dimensional subalgebra.
Since we assume $\g{h}\neq\g{h}_0$, we may write $\g{h}=\operatorname{span}\{A,B\}$ for $A=a_1 E_1+a_3 E_3$ and $B=b_1 E_1+E_2+b_3 E_3$.
Let $g=\Exp(C)$, where $\Exp$ denotes the Lie exponential map and $C=-\frac{b_3}{\lambda_3} E_1+\frac{b_1}{\lambda_1} E_3$.
Then, $\g{h}$ is a subalgebra of $\g{sol}_3$ if and only if $\tilde{\g{h}}=\Ad_g\g{h}=\operatorname{span}\{\Ad_g(A),\Ad_g(B)\}$ is also a subalgebra, where $\Ad_g$ is the adjoint map.
Since $\operatorname{span}\{E_1,E_3\}$ is abelian, we have $\Ad_g(A)=A$.
We compute $[C,B]=-b_3 E_3-b_1E_1$ and hence
\[
	\Ad_g(B)=e^{\ad_{C}}B=B+[C,B]+\frac{1}{2!}[C,[C,B]]+\dots=B+[C,B]=E_2.
\]
Thus, $\g{h}$ is a subalgebra if and only if $\tilde{\g{h}}=\operatorname{span}\{a_1 E_1+a_3 E_3, E_2\}$ is a subalgebra.
Taking now $A=a_1 E_1+a_3 E_3$ and $B=E_2$, condition \eqref{subalgebra condition unimodular} is equivalent to $(-a_3)^2\lambda_1+(a_1)^2\lambda_3=0$.
By rescaling the basis of $\tilde{\g{h}}$ so that $a_1=\sqrt{\lambda_1}$, we get $a_3=\pm\sqrt{-\lambda_3}$, and therefore {$\tilde{\g{h}}=\operatorname{span}\{\sqrt{\lambda_1} E_1\pm \sqrt{-\lambda_3} E_3, E_2\}$}.
Note that the map $\varphi\colon \g{sol}_3\rightarrow\g{sol}_3$ given by
\[
	E_1\longmapsto E_1,\qquad
	E_2\longmapsto -E_2,\qquad
	E_3\longmapsto -E_3,
\]
is an automorphism of $\g{sol}_3$ preserving the metric and maps {$\operatorname{span}\{\sqrt{\lambda_1} E_1+ \sqrt{-\lambda_3} E_3, E_2\}$} to {$\operatorname{span}\{\sqrt{\lambda_1} E_1- \sqrt{-\lambda_3} E_3, E_2\}$}.

Thus, any $2$-dimensional subalgebra of $\g{sol}_3$ is, up to conjugacy and isometric automorphism, either the abelian subalgebra $\g{h}_0=\operatorname{span}\{E_1,E_2\}$ or the non-abelian $\g{h}_1=\operatorname{span}\{\sqrt{\lambda_1} E_1+ \sqrt{-\lambda_3} E_3, E_2\}$.

\subsection{Codimension one subalgebras of non-unimodular Lie algebras}
Consider a non-unimodular Lie algebra $\g{g}=\R^2\rtimes_L\R $  with orthonormal basis $\{E_1,E_2,E_3\}$  such that
\begin{equation*}\label{nonunimodular brackets 2}
	[E_1,E_2]=0,\quad [E_2,E_3]=(1-\alpha)\beta E_1+(\alpha-1)E_2,\quad [E_3,E_1]=(1+\alpha )E_1+(1+\alpha )\beta  E_2.
\end{equation*}
with $\alpha ,\beta \geq0$.
We will assume $\alpha \neq 0, 1$, as otherwise  $\dim(\Isom(G))\geq 4$
(see Remark~\ref{remark isometry group nonunimodular}).

Let $\g{h}$ be a codimension one subalgebra of $\g{g}$ different from the abelian subalgebra $\g{h}_0=\operatorname{span}\{E_1, E_2\}$. Then we can write $\g{h}=\operatorname{span}\{A,B\}$ for $A=a_1 E_1 +a_2 E_2$ and $B=b_1 E_1 +b_2 E_2+ E_3$.
Let $g=\Exp(sE_1+tE_2)$ for $s,t\in \R$. Then, $\g{h}$ is a subalgebra of $\g{g}$ if and only if $\tilde{\g{h}}=\Ad_g\g{h}=\operatorname{span}\{\Ad_g(A),\Ad_g(B)\}$ is a subalgebra of $\g{g}$.
Since $\operatorname{span}\{E_1,E_2\}$ is abelian, we have  $\Ad_g(A)=A$.
We compute 
\[
[sE_1+tE_2,B]=s[E_1, E_3]+t[E_2,E_3]=(-s(1+\alpha )+t\beta (1-\alpha ))E_1+(-s\beta (1+\alpha )+t(\alpha-1 ))E_2,
\]
and hence
\[
\begin{aligned}
	\Ad_g(B)&=e^{\ad_{sE_1+tE_2}}B=B+[sE_1+tE_2,B]+\frac{1}{2!}[sE_1+tE_2,[sE_1+tE_2,B]]+\dots\\
	&=(b_1-s(1+\alpha )+t\beta (1-\alpha ))E_1+(b_2-s\beta (1+\alpha )+t(\alpha-1))E_2+E_3.
\end{aligned}
\]
Note that, since $(1+\beta ^2)(1-\alpha ^2)\neq 0$, there exist unique $s_0,t_0\in \R$ such that
\[
\left\{ 
\begin{array}{l}
	b_1-s_0(1+\alpha )+t_0\beta (1-\alpha )=0, \\
	b_2-s_0\beta (1+\alpha )+t_0(\alpha-1 )=0.
\end{array} 
\right.
\]
Thus, for $g=\Exp(s_0E_1+t_0E_2)$ we get that $\tilde{\g{h}}=\Ad_g\g{h}=\operatorname{span}\{a_1 E_1 +a_2 E_2,E_3\}$.
Now, $\tilde{\g{h}}$ is a subalgebra if and only if
\begin{equation}\label{subalgebra condition nonunimodular}
	0=\langle[a_1 E_1 +a_2E_2,E_3],(a_1 E_1 +a_2E_2)\wedge E_3\rangle=(1+\alpha)\beta {a_1}^2+(1-\alpha)\beta {a_2}^2-2\alpha a_1a_2.
\end{equation}

Note that this is the equation of a degenerate conic passing through the origin, and hence, one can study it according to its discriminant $\alpha^2+\beta^2(\alpha^2-1)=1-\det L$:
\begin{itemize}
	\item If $1-\det L<0$, equation \eqref{subalgebra condition nonunimodular} has only the real solution $a_1=a_2=0$.
	This contradicts the assumption that $\tilde{\g{h}}$ is $2$-dimensional.
	
	\item If $1-\det L>0$, equation \eqref{subalgebra condition nonunimodular} corresponds to a pair of non-coincident intersecting lines, and thus we get two distinct subalgebras.
	For $\beta=0$, these are $\g{h}_1=\operatorname{span}\{E_1,E_3\}$ and $\g{h}_2=\operatorname{span}\{E_2,E_3\}$.
	For $\beta\neq0$, the solutions are given by 
	\[
	a_2=c_\pm a_1\qquad\text{with}\qquad c_\pm=\frac{\alpha\pm\sqrt{1-\det L}}{(1-\alpha)\beta}.
	\]
	The corresponding subalgebras of $\g{g}$ will be denoted by $\g{h}_{\pm}=\operatorname{span}\left\{E_1 +c_{\pm} E_2,E_3\right\}$.
	
	\item If $1-\det L=0$, the only solution to \eqref{subalgebra condition nonunimodular} corresponds to the straight line $\alpha a_1=(1-\alpha)\beta a_2$. We cannot have $\beta=0$, since this contradicts $\det L=1$. Thus $\beta\neq 0$, and we get the subalgebra $\g{h}_{+}=\g{h}_-$ given above, with $c_+=c_-=\alpha/((1-\alpha)\beta)$.
\end{itemize}
Thus, we have proved that any $2$-dimensional Lie subalgebra of $\g{g}$ is conjugate to the abelian subalgebra $\g{h}_0$, to one of the non-abelian subalgebras $\g{h}_1$, $\g{h}_2$ (arising when $\beta=0$ and hence $\det L<1$), to the non-abelian subalgebra $\g{h}_{+}$ (when $\beta\neq 0$ and $\det L=1$), or to $\g{h}_+$, $\g{h}_-$ (if $\beta\neq 0$ and $\det L< 1$). 

Now we will show that no two of these subalgebras are conjugate. Fix one of the non-abelian subalgebras $\g{h}\in\{\g{h}_1,\g{h}_2,\g{h}_+,\g{h}_-\}$. First note that $\g{h}=\mathrm{span}\{X,E_3\}$, for some $X\in\mathrm{span}\{E_1,E_2\}$. Since $\g{h}$ is a subalgebra and $\ad_{E_3}$ leaves $\mathrm{span}\{E_1,E_2\}$ invariant, it follows that $X$ is an eigenvector of $\ad_{E_3}$. Let $Y=sE_1+tE_2+rE_3$, $s,t,r\in\R$, and let $g=\Exp(Y)$ be an arbitrary element in a neighborhood of the identity element of $G$. Then $[Y,X]$ is proportional to $X$. Hence $\Ad_g(X)=e^{\ad_Y}X$ is also proportional to $X$. Since any neighborhood of the identity of $G$ generates $G$, an arbitrary element $g\in G$ can be written as $g=g_1\cdots g_\ell$, with $g_i=\Exp(Y_i)$, $i=1,\dots,\ell$, as before. But then $\Ad_g(\g{h})$ would contain $X$. However $X$ does not belong to any of the subalgebras $\g{h}_1,\g{h}_2,\g{h}_+,\g{h}_-$ different from $\g{h}$ itself. Hence, $\g{h}$ cannot be conjugate to any of the other subalgebras that we have determined. Finally, since we are assuming that the isometry group of $G$ is $3$-dimensional, according to \cite{CoRe} we have $\Isom(G,\langle\cdot,\cdot\rangle)=G$. Thus, by virtue of Proposition~\ref{prop:1-1}, the number of orbit equivalence classes of cohomogeneity one actions on $G$ is one in the case $\det L>1$, two  when $\det L=1$, and three  in the case $\det L<1$.


\section{Cohomogeneity two polar actions}\label{sec:cohom_2}
The aim of this section is to derive the classification of polar actions of cohomogeneity two on simply connected $3$-dimensional homogeneous spaces, up to orbit equivalence. We will first argue that such actions only exist on metric Lie groups with $3$-dimensional isometry group or on symmetric spaces. Then, in Subsections~\ref{subsec:polar_groups} and~\ref{subsec:polar_symmetric} we will discuss each of these two cases separately, thus proving Theorem~\ref{th:c2}.

Let $M$ be a simply connected $3$-dimensional Riemannian homogeneous space, and let $H$ be a connected Lie group acting polarly on $M$ with cohomogeneity two. Thus, there exists a totally geodesic submanifold $\Sigma$ of $M$ that intersects all $H$-orbits and every such intersection is orthogonal. Since the cohomogeneity is two, the section $\Sigma$ is a surface.

The existence of a totally geodesic section $\Sigma$ imposes certain restrictions on $M$. Totally geodesic surfaces in simply connected $3$-dimensional homogeneous spaces have been classified; we refer to the work of Manzano and Souam~\cite{MaSo:mathZ}, where the classification of the more general family of totally umbilical surfaces was achieved. 
It follows from this classification that, if $M$ admits a totally geodesic surface, then it is isometric to a Riemannian symmetric space or to a metric Lie group with a $3$-dimensional isometry group.
In more detail, if $\Sigma\subset M$ is a totally geodesic surface, then we have the following possibilities, see~\cite[Theorems~3.8 and~4.4]{MaSo:mathZ}, where we use the notation of Section~\ref{sec:structure}:

\begin{itemize}
	\item $M$ is isometric to $G=\mathsf{Sol}_3$ with its standard metric, that is, $\lambda_1+\lambda_3=0=\lambda_2$, and $\Sigma$ is an integral surface of one of the distributions spanned by $\{ E_1+E_3,E_2\}$ or $\{E_1- E_3,E_2\}$, all of whose leaves are totally geodesic.
	
	\item $M$ is isometric to $G= \widetilde{\mathsf{SL}_2}(\R)$, $\lambda_3<0<\lambda_2=\lambda_1+\lambda_3<\lambda_1$, and $\Sigma$ is an integral surface of one of the distributions spanned by $\{\sqrt{\lambda_1}\, E_1+\sqrt{-\lambda_3}\,E_3,E_2\}$ or $\{\sqrt{\lambda_1}\, E_1-\sqrt{-\lambda_3}\,E_3,E_2\}$, all of whose leaves are totally geodesic.
	
	\item $M$ is isometric to a non-unimodular Lie group $G=\R^2\rtimes \R$ with $\alpha\neq 0,1$, $\beta=0$, and $\Sigma$ is an integral surface of one of the distributions spanned by $\{ E_1,E_3\}$ or $\{E_2,E_3\}$, all of whose leaves are totally geodesic.
	
	\item $M$ is isometric to a symmetric space, namely $\R^3$, $\mathbb{S}^3$, $\mathbb{H}^3$, $\mathbb{S}^2\times\R$ or $\mathbb{H}^2\times\R$.
\end{itemize}

\begin{remark}\label{rem:only_1_tg}
	Note that each distribution $\g{s}$ in the previous first three items is a subalgebra  of $\g{g}$, and thus its integral surfaces define a foliation of $G$ whose leaves are homogeneous and mutually congruent. In fact, the integral surface through $g\in G$ is $gS$ (namely, the image of the integral surface $S$ under the left translation by $g$), where $S$ is the connected Lie subgroup of $G$ with Lie algebra $\g{s}$. 	It follows from the study of cohomogeneity one actions in the previous section that the choice of sign in the cases $G=\mathsf{Sol}_3$ and $G= \widetilde{\mathsf{SL}_2}(\R)$  gives rise to congruent foliations by mutually congruent totally geodesic surfaces. We warn the reader that the foliations given by integral submanifolds of a Lie subalgebra (as the ones considered in this remark) are in general different from the foliations by orbits of the corresponding group action (as the ones we will consider in Section~\ref{sec:geometry}); the former have  mutually congruent but possibly non-equidistant leaves, whereas the latter have  equidistant but in general non-mutually congruent leaves.
\end{remark}

\subsection{Polar actions on metric Lie groups}\label{subsec:polar_groups}
Suppose that $M$ is isometric to a Lie group $G$ endowed with a left-invariant metric $\langle\cdot,\cdot\rangle$, and such that $\dim(\Isom(G))=3$.
As discussed in Section~\ref{sec:c1}, any isometric action of a connected Lie group on $M$ is equivalent to the action by left translations of some connected subgroup $H\subset G$ on $G$.
Suppose that $H$ acts polarly on $G$ with section $\Sigma$. Without restriction of generality, we can assume that $\Sigma$ passes through the identity element $e$ of $G$, since polar actions admit sections through all points. As argued above, $\Sigma$ is a homogeneous surface, so $\Sigma=S\cdot e=S$ for some connected $2$-dimensional Lie subgroup $S$ of $G$.
Then, since $\Sigma$ intersects $H\cdot e=H$ orthogonally at $e$, $H$ must be the connected Lie subgroup of $G$ with Lie algebra $\g{h}=\g{s}^\perp$, the orthogonal complement of $\g{s}$ in~$\g{g}$. 

The following criterion, due to D\'iaz-Ramos and Kollross, gives us necessary and sufficient conditions for the action of $H$ to be polar.
\begin{lemma}{\cite[Corollary 6]{DRK:dga}}\label{polaritycriterium}
	Let $M$ be a complete connected Riemannian manifold, let $\Sigma$ be a connected totally geodesic embedded submanifold
	of $M$, and let $p\in\Sigma$. 
	Then a proper isometric action of a Lie group $H$ on $M$ is polar with section $\Sigma$ if and only if 
	\begin{enumerate}[\rm(1)]
		\item\label{pcond1} $T_p\Sigma\subset\nu_p (H\cdot p)$, where $\nu_p (H\cdot p)$ is the normal space of the orbit $H\cdot p$ at $p$,
		\item\label{pcond2} the slice representation $H_p\to \mathsf{O}(\nu_p (H\cdot p))$, $h\mapsto h_{*p}\vert_{\nu_p (H\cdot p)}$, of the isotropy subgroup $H_p$ on $\nu_p (H\cdot p)$ is polar with section $T_p\Sigma$, and 
		\item\label{pcond3} $\nabla_v X^*\in \nu_p \Sigma$ for all $v\in T_p\Sigma$ and all Killing fields $X^*$ induced by the $H$-action on~$M$.
	\end{enumerate}
\end{lemma}
We recall that $X^*$ is the fundamental vector field associated with $X\in\g{h}$, which is given by $X^*_p=\frac{d}{dt}\vert_{t=0}\Exp(tX)\cdot p$, for each $p\in M$. If $M$ is a Lie group $G$ with a left-invariant metric, and $X\in\g{g}$, then $X^*$ is precisely the right invariant vector field on $G$ with $X^*_e=X$.

Observe that, in our case, by the relation between $\Sigma$ and $H$, conditions \ref{pcond1} and \ref{pcond2} of Lemma~\ref{polaritycriterium} are automatically satisfied at $p=e$ (note that $H_e$ is trivial).
Thus, the $H$-action is polar if and only if $\langle \nabla_v X^*, w\rangle=0$ for all $X\in\g{h}$ and $v,w\in T_e\Sigma$.
Since $\Sigma$ is homogeneous, we have $T_g\Sigma=\{V^*_g\colon V\in\g{s}\}$ for all $g\in G$, and so we may take $v=V^*_e$ and $w=W^*_e$ for some $V,W\in\g{s}$.
We can now apply Koszul's formula in terms of Killing fields~\cite[Equation~(3.4)]{ORT} to the Killing fields $V^*$, $X^*$, $W^*$ of $(G,\langle\cdot,\cdot\rangle)$:
\[
\begin{aligned}
	2\langle \nabla_{V^*} X^*, W^*\rangle_e&=\langle [V^*,X^*], W^*\rangle_e+\langle [V^*,W^*], X^*\rangle_e+\langle [X^*,W^*], V^*\rangle_e\\
	&=\langle [X,V]^*, W^*\rangle_e+\langle [W,V]^*, X^*\rangle_e+\langle [W,X]^*, V^*\rangle_e\\
	&=\langle [X,V], W\rangle+\langle [W,X], V\rangle,
\end{aligned}
\]
where in the second equality we have used that $Y\in\g{g}\mapsto Y^*\in\g{X}(G)$ is a Lie algebra anti-homomorphism, and in the third one we took into account that $Y^*_e=Y$, for each $Y\in\g{g}$, and, since $V,W\in\g{s}$, we also have $[W,V]^*_e\in\g{s}^*_e=\g{s}$, which is orthogonal to $X^*_e\in\g{h}^*_e=\g{h}$. 

Thus, the action of $H$ is polar if and only if $\langle [X,V], W\rangle+\langle [W,X], V\rangle=0$ for all $X\in \g{h}$ and all $V,W\in \g{s}$.
Note that by linearity and the skew-symmetry of the brackets, and since $\dim \g{h}=1$, it will be enough to check this condition for an ordered basis $\{V,W\}$ of $\g{s}$ and a nonzero $X\in \g{h}$.

We will finally proceed to check the previous condition for the groups that admit totally geodesic surfaces mentioned at the beginning of the section. The groups that do not have totally geodesic surfaces cannot admit cohomogeneity two polar actions.

Let $G=\mathsf{Sol}_3$ or $G= \widetilde{\mathsf{SL}_2}(\R)$, with $\lambda_2=\lambda_1+\lambda_3$. By Remark~\ref{rem:only_1_tg}, there is only one choice for $\g{s}$ up to congruence of the associated totally geodesic foliations. Thus, we have $\g{s}=\operatorname{span}\{V,W\}$ and $\g{h}=\R X$,
where $V=\sqrt{\lambda_1}\, E_1+\sqrt{-\lambda_3}\, E_3$, $W=E_2$ and $X=\sqrt{-\lambda_3}\, E_1-\sqrt{\lambda_1}\, E_3$.
Now,
\[
\langle [X,V], W\rangle+\langle [W,X], V\rangle=2 (\lambda_3^2-\lambda_1^2).
\]
Thus, the action of $H$ is polar if and only if $\lambda_2=\lambda_1+\lambda_3=0$, that is, if $G=\mathsf{Sol}_3$ with its standard metric.
In this case, we have $\g{h}=\R(E_1-E_3)$. This yields the first entry in Table~\ref{table C2:groups}.

For the non-unimodular group $G=\R^2\rtimes \R$, with $\alpha\neq 0,1$, $\beta=0$, we have two possibilities: either $\g{s}=\operatorname{span}\{E_1,E_3\}$ and $\g{h}=\R E_2$ or $\g{s}=\operatorname{span}\{E_2,E_3\}$ and $\g{h}=\R E_1$.
We have
\[
\begin{aligned}
\langle [E_2,E_1], E_3\rangle+\langle [E_3,E_2], E_1\rangle&=\langle(1-\alpha) E_2,E_1\rangle=0, \\
\langle [E_1,E_2], E_3\rangle+\langle [E_3,E_1], E_2\rangle&=\langle(1+\alpha) E_2,E_1\rangle=0.
\end{aligned}
\]
So in both cases the action of $\g{h}$ is polar, and we obtain the second row of Table~\ref{table C2:groups}.

\subsection{Polar actions on symmetric spaces}\label{subsec:polar_symmetric}
Suppose now that $M$ is a $3$-dimensional symmetric space, and write $M=G/K$, where $G$ is, up to a covering, the connected component of the identity of the isometry group of $M$, and $K$ is the stabilizer at some fixed basepoint $o\in M$.
Let $s_o$ be the geodesic reflection of $M$ around $o$.
Then, $\sigma\colon g\mapsto s_o g s_o$ defines an involutive automorphism of $G$, and its differential $\sigma_*$ is an involution of $\g{g}$.
The corresponding decomposition of $\g{g}$ into the $\pm1$-eigenspaces can be written as $\g{g}=\g{k}\oplus\g{p}$, where $\g{k}=\operatorname{Fix}(\sigma_*)$ and $\g{p}$ can be naturally identified with $T_o M$ via $X\in\g{p}\mapsto X^*_o\in T_oM$. We will consider $\g{p}$ endowed with the inner product induced by the metric on $M$ via the isomorphism $\g{p}\cong T_oM$. We also have $[\g{k},\g{p}]\subset{\g{p}}$ and $[\g{p},\g{p}]\subset\g{k}$. See~\cite[Chapter~4, \S3]{Helgason} or \cite[Chapter~6]{Ziller} for further information on the basics of symmetric spaces.

Totally geodesic submanifolds of symmetric spaces can be characterized by their tangent spaces. A subspace $\g{v}$ of $\g{p}$ is called a \emph{Lie triple system} if $[X,[Y,Z]]\in\g{v}$ for all $X,Y,Z\in \g{v}$.
Then, if $\Sigma$ is a totally geodesic submanifold of $M$, its tangent space $T_o\Sigma\subset \g{p}$ is a Lie triple system.
Conversely, if $\g{v}$ is a Lie triple system, then $\exp(\g{v})$ is a totally geodesic
submanifold of $M$, where $\exp$ denotes the Riemannian exponential map; see~\cite[Chapter~4, \S7]{Helgason}.

	\begin{table}[h]
	\centering
	\begin{tabular}{lll}
		\toprule
		$M$ & $H$ & Orbits  \\ \midrule[0.08em]
		\multirow{2}{*}{$\R^3$} & $\R$ & Parallel straigth lines \\ \cline{2-3}
		& $\mathsf{SO}_2$ & An axis of fixed points and circles around it  \\  \hline
		$\mathbb{S}^3(\kappa)$ & $\mathsf{SO}_2$ &  A great circle of fixed points and circles around them\\ \hline
		\multirow{3}{*}{$\H^3(\kappa)$} & $\mathsf{SO}_2$ & A geodesic of fixed points and circles around it \\ \cline{2-3}
		& $A\cong \R$ & A geodesic and equidistant curves to it \\ \cline{2-3}
		& $N'\cong \R$ & Horocycle foliation \\ \hline
		\multirow{2}{*}{$\mathbb{S}^2(\kappa)\times\R$} & $\{\id\}\times\R$ & Vertical lines \\ \cline{2-3}
		& $\mathsf{SO}_2\!\times\!\{0\}$ & Two poles and parallels in each $\mathbb{S}^2(\kappa)\times\{t_0\}$ \\ \hline
		\multirow{3}{*}{$\H^2(\kappa)\times\R$} & $\{\id\}\times\R$ & Vertical lines \\ \cline{2-3}
		& $K\times\{0\}$ & Fixed point and geodesic spheres around it in each $\H^2(\kappa)\times\{t_0\}$ \\ \cline{2-3}
		& $A\times\{0\}$ & A geodesic and equidistant curves to it in each $\H^2(\kappa)\times\{t_0\}$\\ \cline{2-3}
		& $N\times\{0\}$ & Horocycle foliation in each $\H^2(\kappa)\times\{t_0\}$. \\ \bottomrule
	\end{tabular}	
	\bigskip	
	\caption{Cohomogeneity two polar actions on $3$-dimensional symmetric spaces.}
	\label{table C2:symm spaces}
\end{table}

\begin{remark}\label{rem:polar_space_forms}
The classification of cohomogeneity two polar actions on the space forms in Table~\ref{table C2:symm spaces} (namely, $\R^3$, $\mathbb{S}^3(\kappa)$ and $\mathbb{H}^3(\kappa)$) is well known. We briefly review it in this remark.

Dadok~\cite{Dadok} (see also~\cite[Theorem~2.3.17]{BCO}) derived the classification of polar actions on round spheres (or equivalently, of polar orthogonal representations). His theorem states that any such action (of cohomogeneity $k$ on $\mathbb{S}^n$) is orbit equivalent to the isotropy representation of a symmetric space (of dimension $n+1$ and rank $k+1$). Since by classification~\cite[pp.~516--518]{Helgason} the only $4$-dimensional symmetric space of rank $3$ is $\mathbb{S}^2\times\R^2$ (up to quotients and duality, which preserve the isotropy representation), we just get the linear action of $\mathsf{SO}_2$ on $\mathbb{S}^3\subset\R^4$ that fixes a point. 

Regarding polar actions on Euclidean spaces, it can be shown that any such action splits as a product of the action of a polar representation and a translational part (cf.~\cite[Theorem~2.5.1]{BCO}). Combining this fact with Dadok's result, it follows that $\R^3$ only admits two cohomogeneity two actions up to orbit equivalence: the action of a $1$-dimensional translational group, or of a $1$-dimensional rotational group.

Finally, Wu~\cite{Wu} derived the classification of polar actions on hyperbolic spaces $\mathbb{H}^n$, see also~\cite[Theorem~13.5.1]{BCO}. Note that $\mathbb{H}^n=G/K$ with $G=\mathsf{SO}_{1,n}^0$ and $K=\mathsf{SO}_n$. The application of Wu's result yields the three actions by $1$-dimensional Lie subgroups of $\mathsf{SO}_{1,3}^0$ described in Table~\ref{table C2:symm spaces} for $\mathbb{H}^3(\kappa)$. The notation in the table is the standard one used in the Iwasawa decomposition of $G=\mathsf{SO}_{1,3}^0$ as a product $K\times A\times N$, where $K\cong\mathsf{SO}_3$ is the stabilizer of the basepoint $o\in \mathbb{H}^3(\kappa)$, $A$ is the abelian connected subgroup of $G$ whose Lie algebra is any $1$-dimensional subspace of $\g{p}$,  and $N$ is certain (nilpotent, but in this case, abelian) subgroup of $G$ whose action on $\mathbb{H}^3(\kappa)$ produces a foliation by horospheres with a common point at infinity. Thus, the three cohomogeneity two polar actions are given by any $\mathsf{SO}_2$-subgroup of $K\cong\mathsf{SO}_3$, by the group $A\cong\R$ and by any $1$-dimensional subgroup $N'\cong\R$ of $N\cong\R^2$. Again, we refer to \cite[\S13.5.1]{BCO} for details.
\end{remark}

In view of Remark~\ref{rem:polar_space_forms}, in what follows we will suppose that $M$ is one of the product spaces $M=M'\times \R$, where $M'=\mathbb{S}^2(\kappa)$ or $M'=\H^2(\kappa)$.
Write $G=L\times\R$ for the identity connected component of the isometry group of $M$, where $L=\mathsf{SO}_3$ for $M'=\mathbb{S}^2(\kappa)$, and $L=\mathsf{SL}_2(\R)$ for $M'=\H^2(\kappa)$.
In both cases, we can take for example $o=((1,0,0),0)\in M'\times \R$, if we regard $\mathbb{S}^2(\kappa)$ and $\H^2(\kappa)$ as submanifolds of the Euclidean space $\R^3$ and the Minkowski space $\R^{2,1}$, respectively. Then we have $K\cong\mathsf{SO}_2$, which is a Lie subgroup of $L$.
The corresponding decomposition at the Lie algebra level can be written as
\[
\g{g}=\g{k}\oplus\g{p}=\g{k}\oplus\g{p}'\oplus\g{r}=\g{l}\oplus\g{r}
\]
where $\g{r}\cong T_o\R$ and $\g{p}'=\g{p}\ominus\g{r}$ is the orthogonal complement of $\g{r}$ in $\g{p}$ with respect to the induced metric in $\g{p}\cong T_o M$. Note that $\g{p}'$ can be thought of as the tangent space to $M'$.

Suppose that $\g{v}$ is a $2$-dimensional subspace of $\g{p}=\g{p}'\oplus\g{r}$, and let $\g{v}_{\g{p}'}$ denote the orthogonal projection of $\g{v}$ onto $\g{p}'$.
If $\dim \g{v}_{\g{p}'}=1$, then $\g{v}=\operatorname{span}\{X\}\oplus\g{r}$ for some $X\in \g{p}'$, and $\g{v}$ is abelian and therefore a Lie triple system.
Otherwise, it is possible to write {$\g{v}=\operatorname{span}\{X+E_1, Y+E_2\}$}, where $\{X,Y\}$ is an orthonormal basis of $\g{p}'$ and $E_1,E_2\in \g{r}$.
Since $[\g{p},\g{p}]\subset\g{k}$ and $M'$ has nonzero curvature, we have  $[X,Y]\in\g{k}\setminus\{0\}$ (see~\cite[Chapter~4, Theorem~4.2]{Helgason}). Hence, because of the skew-symmetry of $\ad|_{\g{k}}$ we get
\[
\begin{aligned}
	[X+E_1,[X+E_1,Y+E_2]]&=[X,[X,Y]]\in\operatorname{span}\{Y\}\setminus\{0\},\\
	[Y+E_2,[X+E_1,Y+E_2]]&=[Y,[X,Y]]\in\operatorname{span}\{X\}\setminus\{0\}.
\end{aligned}
\]
Thus, $\g{v}$ is a Lie triple system if and only if $E_1=E_2=0$.

Suppose now that $H$ is a connected Lie subgroup of $G$ acting polarly on $M$ with section $\Sigma$ at $o$, and write $\g{s}=T_o\Sigma\subset\g{p}$.
According to the previous discussion, we have that either 
\[
\g{s}=\g{p}'\qquad \text{or} \qquad \g{s}=\operatorname{span}\{X\}\oplus\g{r}.
\]
Via the isomorphism $T_oM\cong\g{p}$, the tangent space of the orbit $H\cdot o$ at $o$ can be identified with the projection $\g{h}_\g{p}$ of $\g{h}$ to $\g{p}$ along $\g{k}$, where $\g{h}$ is the Lie algebra of $H$. Since such tangent space must be orthogonal to $\Sigma$ at $o$ by polarity, we have that $\g{h}_\g{p}$ is orthogonal to $\g{s}$, and hence $\g{h}\subset\g{k}\oplus(\g{p}\ominus\g{s})$, where $\ominus$ denotes orthogonal complement.

If $\g{s}=\g{p}'$, we have $\g{h}\subset\g{k}\oplus\g{r}$.
Suppose that there exists a vector $X\in \g{h}$ with nontrivial projection $T=\pi_{\g{k}}(X)\neq 0$ onto $\g{k}$.
By \eqref{pcond3} in Lemma \ref{polaritycriterium}, we must have $\langle\nabla_{V^*} X^*, W^*\rangle_{o}=0$ for all $V,W\in\g{p}'$.
Now, taking $W=[T,V]\in\g{p}$ and using the fact that $(\nabla_{V^*}X^*)_o=[V^*,X^*]_o=-[V,X]^*_o$ (see~\cite[Proposition~6.34]{Ziller}), we get
\[
\langle\nabla_{V^*} X^*, W^*\rangle_{o}=\langle[X,V]^*, W^*\rangle_{o}=\langle\pi_{\g{p}}([X,V]), W\rangle_{o}=\langle[T,V], W\rangle_{\g{p}}=||W||^2,
\]
which is nonzero if $V\neq 0$.
Thus, we must have $\g{h}=\g{r}$.
It is immediate to see that the corresponding action of $H=\R$, which is given by vertical translations on $M=M'\times\R$, is polar. This yields the first examples of actions on the spaces $\mathbb{S}^2(\kappa)\times\R$ and $\H^2(\kappa)\times \R$ in Table~\ref{table C2:symm spaces}.

Assume now that $\g{s}=\operatorname{span}\{X\}\oplus\g{r}$, for some $X\in\g{p}'$, and so $\g{h}\subset \g{k}\oplus(\g{p}'\ominus X)\subset\g{l}$, where $\g{l}$ is the Lie algebra of $L\in\{\mathsf{SO}_3,\mathsf{SL}_2(\R)\}$.
If $\g{h}$ is two-dimensional, then $\g{l}=\g{sl}_2(\R)$, since $\g{so}_3$ does not admit two-dimensional subalgebras. Moreover, $\g{h}$ is (conjugate to) the parabolic subalgebra described in \S\ref{subalgebras sl2r}. It follows from the theory of symmetric spaces of noncompact type (see for example \cite[\S13.2, pp.~340 and 348]{BCO}) that $H$ acts transitively on $\H^2$, contradicting condition \eqref{pcond1} in Lemma \ref{polaritycriterium}, so $\g{h}$ must be one-dimensional.
Let $H_L$ be the connected subgroup of $L$ with Lie algebra $\g{h}$.
Then, $H=H_L\times\{0\}$, and the orbits of $H$ on $M=M'\times \R$ are of the form $(H_L\cdot p, q)$.
Since $H_L$ is one-dimensional and $L$ acts almost effectively on $M'$, $H_L$ acts with cohomogeneity one on $M'$. On the one hand, since in $\mathbb{S}^2(\kappa)$ there is only one cohomogeneity one action up to orbit equivalence, given by the natural $\mathsf{SO}_2$-action, we obtain the second example corresponding to the space $\mathbb{S}^2(\kappa)\times\R$ in Table~\ref{table C2:symm spaces}. On the other hand, there are three cohomogeneity one actions on the hyperbolic plane $\mathbb{H}^2(\kappa)$, up to orbit equivalence, see~\cite[\S13.5.1, pp.~362--363]{BCO}. These are given by the actions of the isotropy group $K=\mathsf{SO}_2$, of the Lie group $A$ of diagonal matrices in $\mathsf{SL}_2(\R)$ with positive entries in the diagonal, and of the Lie group $N$ of upper triangular matrices in $\mathsf{SL}_2(\R)$ with $1$'s in the diagonal. This yields the remaining three actions in Table~\ref{table C2:symm spaces} corresponding to $\mathbb{H}^2(\kappa)\times \R$.


\section{Geometry of the $2$-dimensional homogeneous foliations}\label{sec:geometry}
If $H$ is a Lie group acting properly by isometries on a Riemannian manifold, its orbits are mutually equidistant submanifolds.
In particular, if $H$ acts with cohomogeneity one and no singular orbits, one may recover the geometry of all orbits by determining the geometry of the orbit through a fixed (but arbitrary) basepoint, and then investigating the geometry of its parallel hypersurfaces.
In this section, we will make use of this idea to study the geometry of homogeneous hypersurfaces of $3$-dimensional metric Lie groups with $3$-dimensional isometry group, thus proving in particular the geometric information in the last column of Table~\ref{table C1}.

Let $G$ be a simply connected $3$-dimensional metric Lie group with $\dim( \Isom(G))=3$. Recall that, by Proposition~\ref{prop:1-1}, there is a one-to-one correspondece between cohomogeneity one actions on $G$ up to orbit equivalence and Lie subalgebras of $\g{g}$ up to conjugacy and isometric automorphisms. In particular, cohomogeneity one actions on $G$ are induced by connected $2$-dimensional subgroups of $G$. 

Thus,  let $H$ be a codimension one Lie subgroup of $G$.
Then, the orbit through the identity element $e$ is precisely $H$. One can compute its shape operator at $e$ with a simple computation at the Lie algebra level.
Note that the orbits of $H$ are precisely the right cosets $H\cdot g$, with $g\in G$. Thus, if $v\in T_g G$ is a normal vector to $H\cdot g$ at a point $g\in G$, then $(L_h)_{*} v$ is a normal vector to $H\cdot g$ at $hg$.
Thus, the left-invariant field $\xi\in\g{g}$ with $\xi_g=v$ is a unit normal field to $H\cdot g$.

Let $\gamma$ be a unit-speed normal geodesic to $H$ with $\gamma(0)=e$ and $\gamma'(0)=v$. Since $G$ is a Riemannian homogeneous space, it is geodesically complete, so we can assume that $\gamma$ is defined in all $\R$.
Let $\xi\in \g{g}$ be the left-invariant field with $\xi_e=v$, and write $H^t$ for the parallel displacement of $H$ in the direction of $\xi$ at distance $t$, that is, $H^t=\{\exp_h(t\xi_h):h\in H\}$ is the parallel (or equidistant) surface to $H$ at distance $t$, where $\exp$ is the Riemannian exponential map.
Then, $H^t=H\cdot \gamma(t)$, and the left-invariant vector field $\xi^t$ with $\xi^t_{\gamma(t)}=\gamma'(t)$ is a unit normal field to $H^t$. In order to calculate the shape operator $S^t$ of $H^t$ with respect to $\xi^t$, we just need to determine the tangent vector $\gamma'$ to the normal geodesic $\gamma$, and then use the formulas \eqref{unimodular connection} and \eqref{eq:LC_non_unimodular} for the Levi-Civita connection of $G$ in terms of left-invariant fields to calculate $S^t=-\nabla_{\cdot\,} \xi^t$ for each $t\in\R$.

Before proceeding with a case-by-case analysis, we will calculate the system of ordinary differential equations defining the geodesic equation on $G$, depending on whether $G$ is unimodular or non-unimodular. Thus, let $\{E_1,E_2,E_3\}$ be some left-invariant orthonormal frame, and write
\[
\gamma'(t)=x(t) E_1 +y(t) E_2 +z(t) E_3,
\]
for some real functions $x$, $y$, $z$.
Denoting by $\nabla$ the Levi-Civita connection of $G$, we have
\begin{equation*}\label{condition geodesics - general case}
0=\nabla_{\gamma'}\gamma'=x' E_1+ x\nabla_{\gamma'}E_1+y' E_2+ y\nabla_{\gamma'}E_2+z' E_3+ z\nabla_{\gamma'}E_3.
\end{equation*}

If $G$ is unimodular and $\{E_1,E_2,E_3\}$ is an orthonormal frame as in \eqref{unimodular brackets}, we can compute using~\eqref{unimodular connection}
\[
\begin{aligned}
	\nabla_{\gamma'}E_1&=x\nabla_{E_1}E_1 +y\nabla_{E_2}E_1+z\nabla_{E_3}E_1=-y\mu_2 E_3+z\mu_3 E_2,\\
	\nabla_{\gamma'}E_2&=x\nabla_{E_1}E_2 +y\nabla_{E_2}E_2+z\nabla_{E_3}E_2=x\mu_1 E_3-z\mu_3 E_1,\\
	\nabla_{\gamma'}E_3&=x\nabla_{E_1}E_3 +y\nabla_{E_2}E_3+z\nabla_{E_3}E_3=-x\mu_1 E_2+y\mu_2 E_1,
\end{aligned}
\]
and so, recalling that $\mu_i=\frac{1}{2}(\lambda_1+\lambda_2+\lambda_3)-\lambda_i$, $i=1,2,3$, we have that $(x,y,z)$ is a solution of the differential equation
\begin{equation}\label{geodesics ODE unimodular}
	\left\{ 
	\begin{array}{l}
		x'+yz(\lambda_3-\lambda_2)=0, \\
		y'+xz(\lambda_1-\lambda_3)=0, \\
		z'+xy(\lambda_2-\lambda_1)=0.
	\end{array} 
	\right.
\end{equation}

For $G$ non-unimodular and $\{E_1,E_2,E_3\}$ satisfying \eqref{brackets nonunimodular}, using~\eqref{eq:LC_non_unimodular}, we get
\[
\begin{aligned}
	\nabla_{\gamma'}E_1&=x\nabla_{E_1}E_1 +y\nabla_{E_2}E_1+z\nabla_{E_3}E_1=z\beta E_2+(x(1+\alpha) +y\alpha\beta )E_3,\\
	\nabla_{\gamma'}E_2&=x\nabla_{E_1}E_2 +y\nabla_{E_2}E_2+z\nabla_{E_3}E_2=-z \beta E_1+(x\alpha\beta +y(1-\alpha)) E_3 \\
	\nabla_{\gamma'}E_3&=x\nabla_{E_1}E_3 +y\nabla_{E_2}E_3+z\nabla_{E_3}E_3=-(x(1+\alpha)+y\alpha\beta)E_1-(x\alpha\beta+y(1-\alpha)) E_2,
\end{aligned}
\]
and hence $(x,y,z)$ is a solution to
\begin{equation}\label{geodesics ODE nonunimodular}
	\left\{ 
	\begin{array}{l}
		x'-(1+\alpha)(x+\beta y)z=0, \\
		y'-(1-\alpha)(y-\beta x)z=0, \\
		z'+(1+\alpha)x^2+2\alpha\beta xy+(1-\alpha)y^2=0.
	\end{array} 
	\right.
\end{equation}

In the rest of this section, we particularize the general framework above to each one of the cohomogeneity one actions classified in Section~\ref{sec:c1} and appearing in Theorem~\ref{th:c1} and Table~\ref{table C1}. This will in particular show the properties in the last column of Table~\ref{table C1}.

\subsection{Homogeneous surfaces of $\widetilde{\mathsf{E}}_2$.}
Recall that the only cohomogeneity one action on $\widetilde{\mathsf{E}}_2=\R^2\rtimes\R$ is the action of the abelian normal subgroup $H=\R^2$.
If $\{E_1,E_2,E_3\}$ is an orthonormal basis of $\g{e}_2$ satisfying the bracket relations described in \eqref{unimodular brackets}, where $\lambda_1>\lambda_2>0=\lambda_3$, then $\g{h}=\operatorname{span}\{E_1,E_2\}$.
It follows from \eqref{unimodular connection} that $\gamma(t)=\Exp(tE_3)$ is a normal geodesic to $H$ through $e$, where $\Exp$ is the Lie exponential map.
Since $\gamma'(t)=E_3$ for all $t\in\R$, the tangent space to $H\cdot \gamma(t)$ at $\gamma(t)$ is the subspace $\g{h}_{\gamma(t)}\subset T_{\gamma(t)}\widetilde{\mathsf{E}}_2$ corresponding to the left-invariant distribution $\g{h}$.
Using~\eqref{unimodular connection}, we have that the shape operator of $H\cdot \gamma(t)$ is given by
\[
	S^t E_1=-\nabla_{E_1} E_3=\mu_1 E_2,\qquad
	S^t E_2=-\nabla_{E_2} E_3=-\mu_2 E_1.
\]
Since $\lambda_3=0$, we have that $\mu_1=-\mu_2=(\lambda_2-\lambda_1)/2\neq 0$, and so the shape operator of $H\cdot \gamma(t)$ at $\gamma(t)$ in terms of the basis $\{E_1,E_2\}$ is given by
\[
S^t\equiv
\begin{pmatrix}
	0&\mu_1\\
	\mu_1&0\\
\end{pmatrix}.
\]
Thus, the principal curvatures of $H\cdot \gamma(t)$ are $\pm\mu_1\neq 0$, with respective principal directions $E_1\pm E_2$. Therefore, all $H$-orbits are minimal, but not totally geodesic.

\subsection{Homogeneous surfaces of $ \widetilde{\mathsf{SL}_2}(\R)$.}
The Lie group $ \widetilde{\mathsf{SL}_2}(\R)$ has a unique codimension one subgroup $H$ up to conjugation.
If $\{E_1,E_2,E_3\}$ is an orthonormal basis of $\g{sl}_2(\R)$ satisfying the bracket relations described in \eqref{unimodular brackets}, where $\lambda_1>\lambda_2>0>\lambda_3$, a representative for the conjugacy class is given by the connected subgroup of $ \widetilde{\mathsf{SL}_2}(\R)$ with Lie algebra $\g{h}=\operatorname{span}\{\sqrt{\lambda_1}\, E_1+\sqrt{-\lambda_3}E_3,E_2\}$.

Let $\gamma$ be a normal geodesic to $H$ with $\gamma(0)=e$ and $\gamma'(0)=\sqrt{\frac{-\lambda_3}{\lambda_1-\lambda_3}}E_1-\sqrt{\frac{\lambda_1}{\lambda_1-\lambda_3}}E_3\in\g{h}^\perp$, and write $\gamma'(t)=x(t) E_1 +y(t) E_2 +z(t) E_3$, $t\in\R$.
Since $(x(t),y(t),z(t))$ is a solution to \eqref{geodesics ODE unimodular}, it satisfies  
$
\frac{xx'}{\lambda_3-\lambda_2}=\frac{yy'}{\lambda_1-\lambda_3}=\frac{zz'}{\lambda_2-\lambda_1}
$.
Integrating, we get that $x^2=\frac{\lambda_3-\lambda_2}{\lambda_1-\lambda_3}y^2+k_1$ and $z^2=\frac{\lambda_2-\lambda_1}{\lambda_1-\lambda_3}y^2+k_2$ for some real constants $k_1, k_2$.
Imposing the initial condition for $\gamma'(0)$ we get
\begin{equation}\label{eq:xz}
	x(t)=\sqrt{\frac{(\lambda_3-\lambda_2)y(t)^2-\lambda_3}{\lambda_1-\lambda_3}},\qquad
	z(t)=-\sqrt{\frac{(\lambda_2-\lambda_1)y(t)^2+\lambda_1}{\lambda_1-\lambda_3}}
\end{equation}
and $y$ satisfies the initial value problem
\[
y'=\sqrt{\left((\lambda_3-\lambda_2)y^2-\lambda_3\right)\left((\lambda_2-\lambda_1)y^2+\lambda_1\right)},\qquad y(0)=0.
\]
In particular, $x(t)>0>z(t)$ for all $t$, so \begin{align*}
V(t)&=\frac{1}{\sqrt{x(t)^2+z(t)^2}} \left(z(t)E_1-x(t) E_3\right),\\
W(t)&=\frac{1}{\sqrt{x(t)^2+z(t)^2}}\left(x(t)y(t) E_1-(x(t)^2+z(t)^2)E_2 +y(t)z(t) E_3\right)
\end{align*}
constitute an orthonormal basis of $T_{\gamma(t)}(H\cdot \gamma(t))$.
Then, the shape operator of $H\cdot \gamma(t)$ is given by
\[
\begin{aligned}
	S^t V(t)=-\nabla_{V(t)}\gamma'(t)&\equiv\frac{xyz(\mu_1-\mu_3)}{x^2+z^2}V-\frac{x^2\mu_3+z^2\mu_1}{x^2+z^2}W,\\
	S^t W(t)=-\nabla_{W(t)}\gamma'(t)&\equiv-\frac{x^2\mu_3+z^2\mu_1}{x^2+z^2}V-\frac{xyz(\mu_1-\mu_3)}{x^2+z^2}W,\\
\end{aligned}
\]
where for the $V$-component of $S^t W$ we used \eqref{eq:xz}, which is equivalent to the self-adjointness of $S^t$. 
Thus, the matrix expression of $S^t$ in terms of the orthonormal basis $\{V(t),W(t)\}$ is
\[
S^t\equiv\frac{1}{x^2+z^2}
\begin{pmatrix}
	xyz(\mu_1-\mu_3)&-(x^2\mu_3+z^2\mu_1)\\
	-(x^2\mu_3+z^2\mu_1)&-xyz(\mu_1-\mu_3)\\
\end{pmatrix}.
\]
Therefore, the orbits of $H$ are all minimal submanifolds of $ \widetilde{\mathsf{SL}_2}(\R)$.
Using~\eqref{eq:xz}, one can verify that the orbit through $e$ is totally geodesic if and only if $\mu_2=0$, or equivalently, $\lambda_2=\lambda_1+\lambda_3$. The other orbits are not totally geodesic. Indeed, since $\lambda_1>\lambda_3$ and $x(t)>0>z(t)$ for all $t$, by~\eqref{geodesics ODE unimodular} we get that $y'(t)>0$ for all $t\in\R$, so $y(t)\neq 0$ for all $t\in\R\setminus\{0\}$ because $y(0)=0$. Then, for $t\neq 0$, $S^t$ does not vanish since $\mu_1-\mu_3=\lambda_3-\lambda_1\neq 0$.

\subsection{Homogeneous surfaces of $\mathsf{Sol}_3$.}

Let  $\{E_1,E_2,E_3\}$ be an orthonormal basis of $\g{sol}_3$ satisfying the bracket relations described in \eqref{unimodular brackets} with $\lambda_1>\lambda_2=0>\lambda_3$. Any $2$-dimensional subalgebra of $\g{sol}_3$ is, up to conjugacy and isometric automorphism, either the abelian subalgebra $\g{h}_0=\operatorname{span}\{E_1,E_3\}$ or the non-abelian subalgebra $\g{h}_1=\operatorname{span}\{\sqrt{\lambda_1} E_1+ \sqrt{-\lambda_3} E_3, E_2\}$, where $\{E_1,E_2,E_3\}$ is an orthonormal basis of $\g{sol}_3$ satisfying the bracket relations in \eqref{unimodular brackets} with $\lambda_1>\lambda_2=0>\lambda_3$.

For $\g{h}_0=\operatorname{span}\{E_1,E_3\}$, $\gamma(t)=\Exp(tE_2)$ is a unit normal geodesic to $H_0$ through $e$, where $H_0$ is the connected subgroup $H_0$ of $\mathsf{Sol}_3$ with Lie algebra $\g{h}_0$.
Thus, $\gamma'(t)=E_2$ for all $t\in\R$, and the tangent space to $H_0\cdot \gamma(t)$ at $\gamma(t)$ is $(\g{h}_0)_{\gamma(t)}$.
We have
\[
	S^t E_1=-\nabla_{E_1} E_2=-\mu_1 E_3, \qquad\qquad
	S^t E_3=-\nabla_{E_3} E_2=\mu_3 E_1.
\]
Since $\lambda_2=0$, then $\mu_1=-\mu_3$, and hence
\[
S^t\equiv
\begin{pmatrix}
	0&\mu_3\\
	\mu_3&0\\
\end{pmatrix}
\]
in terms of the orthonormal basis $\{E_1,E_3\}$ of $T_{\gamma(t)}(H_0\cdot\gamma(t))$.
The principal curvatures of $H_0\cdot\gamma(t)$ are $\pm\mu_3$, and the corresponding principal directions are $E_1\pm E_3$. Hence, all $H_0$-orbits are minimal, but none of them is totally geodesic (since $\mu_3\neq 0$).

Now we consider $\g{h}_1=\operatorname{span}\{\sqrt{\lambda_1}E_1+\sqrt{-\lambda_3}E_3,E_2\}$. Let $\gamma$ denote a unit normal geodesic to the corresponding connected subgroup $H_1$ through $e$ with $\gamma'(0)=\sqrt{\frac{-\lambda_3}{\lambda_1-\lambda_3}} E_1-\sqrt{\frac{\lambda_1}{\lambda_1-\lambda_3}} E_3\in\g{h}^{\perp}_1$. Solving the differential equation \eqref{geodesics ODE unimodular} with initial conditions $x(0)=\sqrt{\frac{-\lambda_3}{\lambda_1-\lambda_3}}$, $y(0)=0$ and $z(0)=-\sqrt{\frac{\lambda_1}{\lambda_1-\lambda_3}}$, we get 
\begin{align*}
\gamma'(t)={}&\sqrt{\frac{-\lambda_3}{\lambda_1-\lambda_3}} \sech\left(t\sqrt{-\lambda_1\lambda_3}\right) E_1+\tanh\left(t\sqrt{-\lambda_1\lambda_3}\right) E_2
\\
&-\sqrt{\frac{\lambda_1}{\lambda_1-\lambda_3}} \sech\left(t\sqrt{-\lambda_1\lambda_3}\right) E_3.
\end{align*}
Let $V(t)=\sqrt{\frac{\lambda_1}{\lambda_1-\lambda_3}} E_1 +\sqrt{\frac{-\lambda_3}{\lambda_1-\lambda_3}} E_3$
and \begin{align*}
W(t)=&{}-\sqrt{\frac{-\lambda_3}{\lambda_1-\lambda_3}}\tanh\left(t\sqrt{-\lambda_1\lambda_3}\right)E_1+\sech\left(t{\sqrt{-\lambda_1\lambda_3}}\right) E_2\\&+\sqrt{\frac{\lambda_1}{\lambda_1-\lambda_3}}\tanh\left(t\sqrt{-\lambda_1\lambda_3}\right)E_3.
\end{align*}
Then, $\{V(t),W(t)\}$ is an orthonormal basis of $T_{\gamma(t)}(H_1\cdot\gamma(t))$, for each $t\in\R$.
We have
\[
\begin{aligned}
	S^t V(t)&=-\nabla_{V(t)} \gamma'(t)=\sqrt{-\lambda_1\lambda_3}\tanh\left(t\sqrt{-\lambda_1\lambda_3} \right)V(t)+\frac{1}{2}\left(\lambda_1+\lambda_3\right) W(t)\\
	S^t W(t)&=-\nabla_{W(t)} \gamma'(t)=\frac{1}{2}\left(\lambda_1+\lambda_3\right) V(t)-\sqrt{-\lambda_1\lambda_3}\tanh\left(t\sqrt{-\lambda_1\lambda_3} \right)W(t).\\
\end{aligned}
\]
Thus, the shape operator of $H_1\cdot\gamma(t)$ with respect to this basis is given by
\[
S^t\equiv
\begin{pmatrix}
	\sqrt{-\lambda_1\lambda_3}\tanh\left(t\sqrt{-\lambda_1\lambda_3} \right)&\frac{1}{2}\left(\lambda_1+\lambda_3\right)\\
	\frac{1}{2}\left(\lambda_1+\lambda_3\right)&-\sqrt{-\lambda_1\lambda_3}\tanh\left(t\sqrt{-\lambda_1\lambda_3} \right)\\
\end{pmatrix}.
\]
The principal curvatures of $H_1\cdot\gamma(t)$ are $\pm\sqrt{\frac{(\lambda_1+\lambda_3)^2}{4}-\lambda_1\lambda_3\tanh^2\left(t\sqrt{-\lambda_1\lambda_3}\right)}$. Note that all $H_1$-orbits are minimal. The only $H_1$-orbit that can be totally geodesic is $H_1\cdot e$, which happens precisely when $\lambda_3=-\lambda_1$.

\subsection{Homogeneous surfaces of non-unimodular groups}
Let $\g{g}$ be a non-unimodular $3$-dimensional Lie algebra satisfying \eqref{brackets nonunimodular} in terms of some orthonormal basis $\{E_1,E_2,E_3\}$.
As in~\S\ref{subsec:unimodular}, we assume $\alpha,\beta\geq 0$, $\alpha\neq 0,1$. Then, we have shown that any $2$-dimensional subalgebra of $\g{g}$ is conjugate to one of the following: the abelian subalgebra $\g{h}_0=\operatorname{span}\{E_1,E_2\}$, the non-abelian subalgebras $\g{h}_1=\operatorname{span}\{E_1,E_3\}$  or $\g{h}_2=\operatorname{span}\{E_2,E_3\}$ (only if $\beta=0$), or the non-abelian subalgebras $\g{h}_{\pm}=\operatorname{span}\{E_1+c_{\pm}E_2,E_3\}$ (only if $\beta\neq 0$).

\subsubsection{The abelian subalgebra $\g{h}_0$}
In this case, we have that $\gamma(t)=\Exp(tE_3)$ is a normal geodesic to $H_0$ through $e$. Thus, $\gamma'(t)=E_3$ for all $t\in\R$, and the tangent space to $H_0\cdot \gamma(t)$ at $\gamma(t)$ is $(\g{h}_0)_{\gamma(t)}$. Hence, using~\eqref{eq:LC_non_unimodular}, we get
\[
S^t E_1=-\nabla_{E_1} E_3=(1+\alpha) E_1+\alpha\beta E_2,\qquad
S^t E_2=-\nabla_{E_2} E_3=\alpha\beta E_1 + (1-\alpha) E_2.
\]
Therefore, the principal curvatures of $H_0\cdot\gamma(t)$ are $1\pm\alpha\sqrt{1+\beta^2}$. In particular, the mean curvature is $1$, and none of the orbits of $H_0$ is minimal.

\subsubsection{The non-abelian subalgebras $\g{h}_1$ and $\g{h}_2$ in the case $\beta=0$} 

For $\g{h}_1=\mathrm{span}\{E_1,E_3\}$, let $\gamma$ be the unit normal geodesic to the corresponding connected Lie subgroup $H_1$ of $G$ with $\gamma'(0)=E_2\in\g{h}_1^\perp$. Solving equation~\eqref{geodesics ODE nonunimodular} we have that $\gamma'(t)=\sech\left((1-\alpha)t\right) E_2-\tanh\left((1-\alpha)t\right) E_3$.
Thus, an orthonormal basis of $T_{\gamma(t)}(H_1\cdot\gamma(t))$ is given by $V(t)=E_1$ and $W(t)=\tanh\left((1-\alpha)t\right) E_2 +\sech\left((1-\alpha)t\right) E_3$. Using~\eqref{eq:LC_non_unimodular}, we get
\[
\begin{aligned}
	S^t V(t)&=-\nabla_{V(t)} \gamma'(t)&=&-(1+\alpha)\tanh\left((1-\alpha)t\right)V(t),\\
	S^t W(t)&=-\nabla_{W(t)} \gamma'(t)&=&-(1-\alpha)\tanh\left((1-\alpha)t\right)W(t),\\
\end{aligned}
\]
so $V(t)$ and $W(t)$ are principal directions of  $H_1\cdot\gamma(t)$. Thus, the $H_1$-orbit through $\gamma(t)$ has constant mean curvature $-\tanh\left((1-\alpha)t\right)\in(-1,1)$. In particular, the only minimal orbit is the one through $e$, which is actually totally geodesic.

For $\g{h}_2=\mathrm{span}\{E_2,E_3\}$, if $\gamma$ is the normal geodesic to $H_2$ with $\gamma'(0)=E_1\in\g{h}_2^\perp$, we have that
$\gamma'(t)=\sech\left((1+\alpha)t\right) E_1-\tanh\left((1+\alpha)t\right) E_3$, and so $V=E_2$ and $W=\tanh\left((1+\alpha)t\right) E_1+\sech\left((1+\alpha)t\right) E_3$
provide an orthonormal basis for $T_{\gamma(t)}(H_2\cdot\gamma(t))$.
In this case,
\[
\begin{aligned}
	S^t V(t)&=-\nabla_{V(t)} \gamma'(t) &=& -(1-\alpha)\tanh\left((1+\alpha)t\right)V(t),\\
	S^t W(t)&=-\nabla_{W(t)} \gamma'(t) &=& -(1+\alpha)\tanh\left((1+\alpha)t\right)W(t).\\
\end{aligned}
\]
Therefore, the orbit $H_2\cdot\gamma(t)$ has constant mean curvature $-\tanh\left((1+\alpha)t\right)\in(-1,1)$, and $H_2\cdot e$ is totally geodesic.

\subsubsection{The non-abelian subalgebras $\g{h}_+$ and $\g{h}_-$ in the case $\beta\neq0$}
Finally, we consider a non-unimodular Lie algebra $\g{g}$ with $\beta\neq 0$, and the Lie subalgebras $\g{h}_{\pm}=\operatorname{span}\{E_1+c_{\pm}E_2,E_3\}$, where  $c_\pm=\frac{\alpha\pm\sqrt{1-(1-\alpha^2)(1+\beta^2)}}{(1-\alpha)\beta}=\frac{\alpha\pm\sqrt{1-\det L}}{(1-\alpha)\beta}$. We will treat both cases $\g{h}_+$ and $\g{h}_-$ simultaneously.
We have that $\frac{1}{\sqrt{1+c_{\pm}^2}} (c_{\pm}E_1- E_2)$ is a unit normal vector to $\g{h}_{\pm}$.
Solving the initial value problem~\eqref{geodesics ODE nonunimodular} with initial conditions $x(0)=\frac{c_{\pm}}{\sqrt{1+c_{\pm}^2}}$, $y(0)=\frac{-1}{\sqrt{1+c_{\pm}^2}}$ and $z(0)=0$ yields
\[
\gamma'(t)=\frac{c_{\pm}}{\sqrt{1+c_{\pm}^2}}\sech\left(\omega t\right) E_1-\frac{1}{\sqrt{1+c_{\pm}^2}}\sech\left(\omega t\right) E_2-\tanh\left(\omega t\right) E_3,
\]
where $\omega =(1-\alpha)(1+\beta c_{\pm})$ is a nonzero constant, since $\alpha\neq 1$ and $1+\beta c_\pm\neq 0$. For each $t\in\R$ we define
\begin{align*}
	V(t)&=\frac{1}{\sqrt{1+c_{\pm}^2}}(E_1+c_{\pm}E_2),
	\\
	W(t)&=\frac{1}{\sqrt{1+c_{\pm}^2}}\tanh\left(\omega t\right) (-c_\pm E_1+E_2)-\sech\left(\omega t\right) E_3,
\end{align*}
so that $\{V(t),W(t)\}$ is an orthonormal basis of $H_{\pm}\cdot\gamma(t)$.
Then we have 
\[
\begin{aligned}
	S^t V(t)&=\bigl(\pm\sqrt{1-\det L}-1\bigr)\tanh(\omega t)V(t) +\beta W(t),\\
	S^t W(t)&=\beta V(t)- \bigl(\pm\sqrt{1-\det L}+1\bigr)\tanh(\omega t) W(t),\\
\end{aligned}
\]
where we have used the relation between $c_\pm$, $\alpha$ and $\beta$ to simplify expressions. The principal curvatures of $H_\pm\cdot \gamma(t)$ are $-\tanh(\omega t)\pm\sqrt{(1-\det L)\tanh^2(\omega t)+\beta^2}$, and its mean curvature  is $-\tanh(\omega t)\in(-1,1)$. Since $\omega \neq 0$, we deduce that the only minimal $H_\pm$-orbit is the one through the identity element, but it is not totally geodesic as $\beta\neq 0$.

\enlargethispage{2\baselineskip}

\begin{thebibliography}{99}
\bibitem[AR04]{AR:acta} U.~Abresch, H.~Rosenberg: A Hopf differential for constant
mean curvature surfaces in $S^2\times\R$ and $H^2\times \R$. \emph{Acta Math.} \textbf{193} (2004), 141--174.	
	
\bibitem[BCO16]{BCO}	J.~Berndt, S.~Console, C.~Olmos: \emph{Submanifolds and holonomy}. Second edition. Monographs and Research Notes in Mathematics, CRC Press, Boca Raton, FL, 2016. 
	
\bibitem[CR15]{CR:book} T.~E.~Cecil, P.~J.~Ryan: \emph{Geometry of hypersurfaces}. Springer Monographs in Mathematics, Springer, 2015.
	
\bibitem[Ch41]{Che} C.~Chevalley: On the topological structure of solvable Lie groups. \textit{Ann.\ Math.}\ \textbf{42} (1941), 668--675.

\bibitem[CR22]{CoRe} A.~Cosgaya, S.~Reggiani: Isometry groups of three-dimensional Lie groups. \textit{Ann.\ Glob.\ Anal.\ Geom.}\ \textbf{61} (2022), 831--845.

\bibitem[Da85]{Dadok} J.~Dadok: Polar coordinates induced by actions of compact Lie groups. \textit{Trans.\ Amer.\ Math.\ Soc.}\ \textbf{288} (1985), 125--137.

\bibitem[DK11]{DRK:dga} J.~C.~D\'iaz-Ramos, A.~Kollross: Polar actions with a fixed point. \textit{Diff.\ Geom.\ Appl.}\ \textbf{29} (2011), 20--25.

\bibitem[DM21]{DVMa:annali} M.~Dom\'inguez-V\'azquez, J.~M.~Manzano: Isoparametric surfaces in $\mathbb{E}(\kappa,\tau)$-spaces. \textit{Ann.\ Sc.\ Norm.\ Super.\ Pisa Cl.\ Sci.\ (5)} \textbf{22} (2021), 269--285.

\bibitem[FGT17]{FGT} F.~Fang, K.~Grove, G.~Thorbergsson: Tits geometry and positive curvature. \emph{Acta Math.}\ \textbf{218} (2017), 1--53.


\bibitem[FH17]{FH} L.~Foscolo, M.~Haskins: New $\mathsf{G}_2$-holonomy cones and exotic nearly K\"ahler structures on $S^6$ and $S^3 \times S^3$. \emph{Ann.\ of Math.\ (2)} \textbf{185} (2017), no.~1, 59--130.

\bibitem[GKR25]{GKRV} C.~Gorodski, A.~Kollross, A.~Rodríguez-Vázquez: Totally geodesic submanifolds and polar actions on Stiefel manifolds. \emph{J.\ Geom.\ Anal.}\ \textbf{35} (2025), article no.~41.

\bibitem[HL12]{HaLee} K.~Y.~Ha, J.~B.~Lee: The isometry groups of simply connected 3-dimensional unimodular Lie groups. \textit{J.\ Geom.\ Phys.}\ \textbf{62} (2012), 189--203.

\bibitem[He78]{Helgason} S.~Helgason: \emph{Differential geometry, Lie groups, and symmetric spaces}. Pure and
Applied Mathematics, 80, Academic Press, Inc., New York-London, 1978.

\bibitem[HH89]{HH:inventiones} W.-T.~Hsiang, W.-y.~Hsiang: On the uniqueness of isoperimetric solutions and imbedded soap bubbles in non-compact symmetric spaces, I. \emph{Invent.\ Math.}\ \textbf{98} (1989), 39--58. 

\bibitem[Ma45]{Mal} A.~Malcev: On the theory of the Lie groups in the large. \textit{Rec.\ Math.\ [Mat.\ Sbornik] N.S.}\ \textbf{16(58)} (1945), no.~2, 163--190.

\bibitem[MS15]{MaSo:mathZ} J.~M.~Manzano, R.~Souam: The classification of totally umbilical surfaces in homogenous $3$-manifolds. \textit{Math.\ Z.}\ \textbf{279} (2015), 257--576.

\bibitem[MP12]{MePe:cont} W.~H.~Meeks, J.~P\'erez: Constant mean curvature surfaces in metric Lie groups. \textit{Geometric Analysis: Partial Differential Equations and Surfaces}, Contemporary Mathematics (AMS)
vol.~570 (2012), 25--110.

\bibitem[Mi08]{Michor} P.~W.~Michor: \emph{Topics in Differential Geometry}. Graduate Studies in Mathematics, vol.~93, American Mathematical Society, 2008.

\bibitem[Mi76]{Mi:adv} J.~Milnor: Curvatures of left invariant metrics on Lie groups. \textit{Adv.\ Math.}\ \textbf{21} (1976), 293--329.

\bibitem[ORT14]{ORT} C.~Olmos, S.~Reggiani, H.~Tamaru: The index of symmetry of compact naturally reductive
spaces. \emph{Math.\ Z.}\ \textbf{277} (2014), 611--628.

\bibitem[Ot24]{Otero} T.~Otero Casal: \emph{Homogeneous hypersurfaces in symmetric spaces}. PhD Thesis, Universidade de Santiago de Compostela.

\bibitem[PT87]{PT} R.~S.~Palais, C.-L.~Terng: A general theory of canonical forms. \emph{Trans.\ Amer.\
Math.\ Soc.}\ \textbf{300} (1987), no.~2, 771--789.

\bibitem[SS25]{SoSa} I.~Solonenko, V.~Sanmart\'in-L\'opez: Classification of cohomogeneity-one actions on symmetric spaces of noncompact type. arXiv:2501.05553.

%

\bibitem[Wu92]{Wu} B.~Wu: Isoparametric submanifolds of hyperbolic spaces. \emph{Trans.\ Amer.\ Math.\ Soc.}\ \textbf{331} (1992), no.~2, 609–626.

\bibitem[Zi10]{Ziller} W.~Ziller: \emph{Lie groups,  representation theory and symmetric spaces}. Lecture notes. Available on \url{https://www2.math.upenn.edu/~wziller/math650/LieGroupsReps.pdf} (latest access on 11 February 2025).

\end{thebibliography}
\end{document}